\documentclass[aos]{imsart}

\RequirePackage{amsthm,amsmath,amsfonts,amssymb}
\RequirePackage[numbers]{natbib}
\RequirePackage[colorlinks,citecolor=blue,urlcolor=blue]{hyperref}
\RequirePackage{graphicx}
\RequirePackage{subcaption} 

\startlocaldefs
\theoremstyle{plain}

\newtheorem{theorem}{Theorem}[section]
\newtheorem{lemma}[theorem]{Lemma}
\theoremstyle{remark}

\newtheorem*{remark}{Remark}

\newtheorem{proposition}{Proposition}
\newcommand{\sgn}{\text{sgn}}
\usepackage{comment}
\endlocaldefs

\usepackage{amsthm}

\newcommand{\diag}{\text{diag}}
\newcommand{\tr}{\text{tr}}

\begin{document}
\begin{frontmatter}
\title{Approximate Risk Minimization over Shrinkage–Thresholding Rules in Normal Mean Estimation}
\runtitle{Approximate Risk Minimization in Normal Mean Estimation}

\begin{aug}
\author[A, B]{\fnms{Wei}~\snm{Jiang}\ead[label=e1]{wei.jiang@uta.edu}\orcid{0000-0001-6120-5278}}
\address[A]{Department of Mathematics,
University of Texas at Arlington, Arlington, TX, US, 76019}
\address[B]{Division of Data Science,
University of Texas at Arlington, Arlington, TX, US, 76019\printead[presep={,\ }]{e1}}

\end{aug}

\begin{abstract}
We develop an approximate risk minimization framework for shrinkage-thresholding estimation in normal mean problems. In the canonical multivariate normal mean model, we introduce a general functional class of estimators that contains classical shrinkage and thresholding behavior, including James–Stein-type and lasso-type rules. We express quadratic risk as a functional over this class, derive optimality conditions for both oracle risk and data-driven approximate risk minimization, and construct a feasible approximate risk criterion from the observed data when the oracle risk is unavailable. The resulting estimator, NOMAD, is obtained by minimizing this approximate risk over the proposed class.

For the canonical model, we develop an approximate risk minimization theory that includes optimizer characterization, sieve-based consistency under regularity conditions, and approximate-risk inequalities relative to benchmark procedures in the admissible class. We then extend the framework to multivariate normal mean estimation with correlated observations, develop both MLE-based and conditional MLE-based constructions, and establish consistency results under regularity conditions. We further apply the framework to linear regression and derive an equivalent penalized regression representation in which the shrinkage-thresholding map induces a data-adaptive penalty, recovering ridge-type and lasso-type behavior as special cases or limiting forms. The results provide a unified risk-based framework for shrinkage, thresholding, and regularization across canonical and correlated normal mean estimation and linear regression.

\end{abstract}



\end{frontmatter}


\section{Introduction}

Shrinkage and thresholding are two central mechanisms for improving high-dimensional estimation under squared error loss and they arose from a long trajectory in frequentist estimation. After Fisher developed the foundations of maximum likelihood estimation between 1912 and 1922, the maximum likelihood estimator (MLE) became a default method because it is consistent and it achieves classical efficiency in regular parametric settings \cite{aldrich1997ra}. In high dimensions, however, Stein showed that the MLE is inadmissible for multivariate normal mean estimation under quadratic loss, and James and Stein exhibited an estimator that uniformly improves the MLE by shrinking the observed vector toward the origin \cite{stein1956inadmissibility,james1961estimation}. This shrinkage principle motivated a broad family of regularized estimators, including ridge-type procedures that shrink continuously toward zero \cite{hoerl1970ridge}. A complementary development arose from sparse modeling, where the lasso penalization shrinks coefficients and it also induces exact zeros through the nonsmooth geometry of the $\ell_1$ norm \cite{tibshirani1996regression}. The resulting thresholding effect can yield lower risk than purely shrinkage-based rules in sparse regimes \cite{hansen2016risk}. Despite these connections, shrinkage and thresholding are often developed through separate methodological lenses, which can obscure their shared risk structure and it can hinder a unified theory that explains how procedures adapt across dense and sparse signal settings.

This paper develops an approximate risk minimization framework for shrinkage and thresholding estimation in normal mean problems. We begin with the canonical multivariate normal mean model and introduce a general functional class of estimators that combines shrinkage and thresholding effects through two functional components. We refer to this class as the general shrinkage-thresholding estimator family (GEST) and we show that it contains classical estimators, such as the maximum likelihood estimator, James--Stein type shrinkage rules, and lasso-type thresholding rules, as special cases. We express quadratic risk as a functional over this class and derive optimality conditions for risk minimization. Oracle risk minimization is infeasible because the risk depends on the unknown mean vector, so we construct a data-dependent approximate risk criterion from the observed data. The resulting estimator is obtained by minimizing this approximate risk within the proposed class and we refer to it as \emph{NOMAD}, which is short for \emph{NOnparametric Minimizer of an Approximate Data-driven Risk criterion}. Throughout the paper we use the shorter term \emph{approximate risk} to denote this data-dependent criterion and we refer to the associated results as \emph{approximate risk minimization theory}.

The approximate risk minimization theory developed for the canonical model characterizes the optimizer within the proposed function class, establishes sieve-based consistency under regularity conditions, and yields approximate-risk inequalities relative to benchmark procedures in the admissible class. These results clarify how the selected rule adapts between shrinkage-dominant and thresholding-dominant regimes and provide a principled route to estimators that interpolate between dense and sparse settings within a unified framework. The canonical normal mean formulation also leads directly to wavelet denoising. Under standard Gaussian noise models, an orthogonal wavelet transform reduces wavelet denoising to normal mean estimation in the coefficient domain \cite{donoho1994ideal,donoho1995adapting}. This connection allows NOMAD to be applied as a coefficient-domain shrinkage-thresholding rule chosen by approximate risk minimization and provides an interpretable illustration of the framework in a classical nonparametric setting.

Many problems involve dependence that precludes reduction to the canonical model. We therefore extend the framework to multivariate normal mean estimation with correlated observations. In this setting, we develop shrinkage-thresholding constructions based on both the maximum likelihood estimator and a conditional maximum likelihood estimator, together with their corresponding data-driven criteria. For the correlated model, we establish consistency results under regularity conditions and clarify the role of the conditional construction as a structurally motivated alternative that incorporates local dependence adjustment before thresholding. We then connect the correlated normal mean formulation to linear regression and show that the resulting estimator admits an equivalent penalized regression representation with a data-adaptive penalty induced by the shrinkage-thresholding map. This representation recovers ridge-type and lasso-type behavior as special cases or limiting forms and provides a unified risk-based perspective on shrinkage, thresholding, and regularization across canonical normal means, correlated normal means, and linear regression.

\paragraph*{Contributions and Organization}
The main contributions and corresponding sections of this work are as follows.
First, we introduce a general functional class of shrinkage-thresholding estimators in the canonical normal mean model that unifies classical shrinkage and thresholding behavior within a single representation (Section~\ref{sec:canonical_model}).
Second, we construct a data-dependent approximate risk criterion (Section~\ref{sec:approx_risk}), define NOMAD as its minimizer over the proposed class (Section~\ref{sec:nomad}), and study its connections with classical estimators (Section~\ref{sec:connections}).
Third, we develop approximate risk minimization theory in the canonical model, including optimizer characterization, sieve-based consistency, and approximate-risk inequalities relative to benchmark procedures in the admissible class (Section~\ref{sec:canonical_theory}).
Fourth, we present wavelet denoising as a direct application of the canonical framework (Section~\ref{sec:wavelets}).
Fifth, we extend the framework to correlated multivariate normal mean estimation via MLE-based and conditional MLE-based constructions and establish consistency results under regularity conditions (Section~\ref{sec:correlated_extension}).
Sixth, we derive a linear regression formulation and an equivalent penalized regression representation with a data-adaptive penalty induced by the shrinkage-thresholding map, recovering ridge-type and lasso-type behavior as special cases or limiting forms (Section~\ref{sec:regression}).
Finally, Section~\ref{sec:simulations} reports simulation studies and Section~\ref{sec:discussion} concludes with a discussion of implications and extensions.

\section{Canonical Normal Mean Model and the GEST Class}\label{sec:canonical_model}

\subsection{Canonical multivariate normal mean model}\label{sec:canonical_model_setup}
We begin with the canonical problem of estimating the mean vector in a multivariate normal distribution,
\begin{equation}
\pmb{Z}\sim N(\pmb{\theta}, \pmb{I}_d),\label{eq:canonical_model}
\end{equation}
where $\pmb{Z}=(Z_1,\ldots,Z_d)^\top$ is the observed vector, $\pmb{\theta}=(\theta_1,\ldots,\theta_d)^\top$ is the unknown mean vector, and $\pmb{I}_d$ is the $d$-dimensional identity matrix. The coordinates of $\pmb{Z}$ are independent under \eqref{eq:canonical_model}. We evaluate an estimator $\hat{\pmb{\theta}}=\hat{\pmb{\theta}}(\pmb{Z})$ under 
quadratic risk
\begin{equation}
L(\hat{\pmb{\theta}},\pmb{\theta})=\|\hat{\pmb{\theta}}-\pmb{\theta}\|_2^2,
\qquad
R(\hat{\pmb{\theta}},\pmb{\theta})=\mathbb{E}_{\pmb{\theta}}\!\left[\|\hat{\pmb{\theta}}-\pmb{\theta}\|_2^2\right]
=\sum_{i=1}^d \mathbb{E}_{\pmb{\theta}}\!\left(\hat{\theta}_i-\theta_i\right)^2.\label{eq:quadratic_risk}
\end{equation}
Our goal is to develop a general class of shrinkage-thresholding estimators for \eqref{eq:canonical_model} and to later select an estimator within this class by approximate risk minimization.

\subsection{General family of estimators with shrinkage and thresholding effects (GEST)}\label{sec:gest}
We first recall several commonly used estimators for \eqref{eq:canonical_model}. The maximum likelihood estimator is $\hat{\pmb{\theta}}_{\mathrm{MLE}}=\pmb{Z}$. The James--Stein estimator takes the form
\begin{equation}
\hat{\pmb{\theta}}_{\mathrm{JS}}
=\left(1-\frac{d-2}{\|\pmb{Z}\|_2^2}\right)\pmb{Z},\label{eq:js}
\end{equation}
and its positive-part modification replaces the shrinkage factor by its positive part $\hat{\pmb{\theta}}_{\mathrm{JS}}^{+}$,
which is known to dominate the original James--Stein rule under quadratic loss \cite{baranchik1964multiple}. Regularization-based estimators provide additional examples. Ridge regression in the canonical model yields the linear shrinkage rule
\begin{equation}
\hat{\pmb{\theta}}_{\mathrm{ridge}}
=\arg\min_{\hat{\pmb{\theta}}}\ \|\pmb{Z}-\hat{\pmb{\theta}}\|_2^2+\lambda\|\hat{\pmb{\theta}}\|_2^2
=\left(1-\frac{\lambda}{1+\lambda}\right)\pmb{Z},\label{eq:ridge}
\end{equation}
and the lasso yields the coordinatewise soft-thresholding rule
\begin{equation}
\hat{\pmb{\theta}}_{\mathrm{lasso}}
=\arg\min_{\hat{\pmb{\theta}}}\ \frac12\|\pmb{Z}-\hat{\pmb{\theta}}\|_2^2+\lambda\|\hat{\pmb{\theta}}\|_1
=\left(1-\frac{\lambda}{|\pmb{Z}|}\right)_+\circ \pmb{Z},\label{eq:lasso}
\end{equation}
where $\lambda>0$, $|\pmb{Z}|=(|Z_1|,\ldots,|Z_d|)^\top$, and $\circ$ denotes the Hadamard product. Ridge and James--Stein shrink continuously toward the origin, while the lasso shrinks and also thresholds small coordinates to exact zero.

Motivated by these forms, we consider a general family of estimators that combines shrinkage and thresholding through two nonnegative functional components. We define the truncated general shrinkage-thresholding estimator (GEST) class by
\begin{equation}
\hat{\pmb{\theta}}_{\mathrm{GEST}}^{+}
=\left(1-\frac{t(|\pmb{Z}|)}{|\pmb{Z}|^{c(|\pmb{Z}|)}}\right)_+\circ \pmb{Z},
\qquad
t(|\pmb{Z}|)\ge 0,\quad c(|\pmb{Z}|)\ge 0.\label{eq:gest_plus}
\end{equation}
Here $t(\cdot)$ is a nonnegative threshold function and $c(\cdot)$ is a nonnegative rate function, and both may depend on the full vector $|\pmb{Z}|$. The shrinkage factor in \eqref{eq:gest_plus} decreases as $t(|\pmb{Z}|)/|\pmb{Z}|^{c(|\pmb{Z}|)}$ increases, and the truncation $(\cdot)_+$ ensures that coordinates with $|Z_i|^{c(|\pmb{Z}|)}\le t(|\pmb{Z}|)$ are set to zero. The rate function $c(\cdot)$ controls how shrinkage decays with the magnitude of a coordinate, while the threshold function $t(\cdot)$ determines which coordinates are removed by thresholding after the scale adjustment $|Z_i|^{c(|\pmb{Z}|)}$.

The class \eqref{eq:gest_plus} includes the estimators above as special cases through particular choices of $(t,c)$. In particular,
\begin{eqnarray}
\mathrm{MLE}&:& t(|\pmb{Z}|)=0,\quad c(|\pmb{Z}|)\ge 0,\nonumber\\
\mathrm{JS}^{+}&:& t(|\pmb{Z}|)=(d-2)/\|\pmb{Z}\|_2^2,\quad c(|\pmb{Z}|)=0,\label{eq:gest_special_cases}\\
\mathrm{ridge}&:& t(|\pmb{Z}|)=\lambda/(1+\lambda),\quad c(|\pmb{Z}|)=0,\nonumber\\
\mathrm{lasso}&:& t(|\pmb{Z}|)=\lambda,\quad c(|\pmb{Z}|)=1.\nonumber
\end{eqnarray}
These examples motivate viewing $t(\cdot)$ and $c(\cdot)$ as functional parameters that encode thresholding and shrinkage behavior within one representation.

\section{Risk and Approximate Risk in the Canonical Model}\label{sec:approx_risk}

This section studies the quadratic risk of the untruncated GEST family in the canonical normal mean model \eqref{eq:canonical_model}. The oracle risk depends on the unknown mean vector and it is not computable from one realization. We therefore construct a data-dependent approximation by combining a risk representation based on Stein's identity with a nonparametric plug-in approximation of $\pmb{\theta}$ based on Tweedie's formula \cite{efron2011tweedie}.

\subsection{Oracle quadratic risk as a functional over GEST}\label{sec:oracle_risk_gest}

\begin{proposition}\label{prop:risk}
The risk of $\hat{\pmb{\theta}}_{\text{GEST}}$ can be written as
\[
R(\hat{\pmb{\theta}}_{\text{GEST}}, \pmb{\theta})
=
d+\int_{\mathbb{R}^d}
H\!\left(\pmb{z}, t, c, \frac{\partial t}{\partial |\pmb{z}|}, \frac{\partial c}{\partial |\pmb{z}|}\right)\,d\pmb{z},
\]
where $\pmb{z}=(z_1,\ldots,z_d)^\top$ and
\begin{eqnarray*}
H\!\left(\pmb{z}, t, c, \frac{\partial t}{\partial |\pmb{z}|}, \frac{\partial c}{\partial |\pmb{z}|}\right)
&=&
-\sum_{i=1}^d \frac{1}{|z_i|^{c(|\pmb{z}|)}}\Bigl(
2|z_i|\frac{\partial t(|\pmb{z}|)}{\partial |z_i|}
-
2t(|\pmb{z}|)|z_i|\ln|z_i|\frac{\partial c(|\pmb{z}|)}{\partial |z_i|}
\\
&&
+
2t(|\pmb{z}|)(1-c(|\pmb{z}|))
-
t(|\pmb{z}|)^2|z_i|^{2-c(|\pmb{z}|)}
\Bigr)
\prod_{j=1}^d \phi(z_j-\theta_j),
\end{eqnarray*}
and $\phi(\cdot)$ is the standard normal density.
\end{proposition}

The proof of Proposition~\ref{prop:risk} is based on Stein's Lemma \cite{stein1981estimation} and is demonstrated in the Supplement \cite{jiang2026nomadS}. Proposition~\ref{prop:risk} shows that the risk depends on $\pmb{\theta}$ through the likelihood factor $\prod_{j=1}^d \phi(z_j-\theta_j)$. This dependence makes direct risk minimization infeasible in practice.

\subsection{Data-dependent approximate risk criterion}\label{sec:approx_risk_criterion}

We construct a data-dependent approximation by replacing the unknown $\pmb{\theta}$ in Proposition~\ref{prop:risk} by a nonparametric plug-in estimate based on Tweedie's formula. This construction uses only the empirical distribution of the observed coordinates and it does not require specification of a parametric prior \cite{efron2011tweedie}.

\begin{lemma}[Robbins (1956) \cite{robbins1992empirical}]\label{lem:tweedie}
Suppose $\theta_1,\ldots,\theta_d$ are sampled from a prior distribution $g(\theta)$ and $\check{z}_i$ is an observation from $N(\theta_i,1)$. Let $f(\check{z}_i)=\int \phi(\check{z}_i-\theta)g(\theta)\,d\theta$ denote the marginal density. Then
\[
\mathbb{E}(\theta_i\mid \check{z}_i)=\check{z}_i+l'(\check{z}_i),
\qquad
l'(\check{z}_i)=\frac{d}{d\check{z}_i}\log f(\check{z}_i).
\]
\end{lemma}

Lemma~\ref{lem:tweedie} motivates approximating $\theta_i$ by $\check{z}_i+l'(\check{z}_i)$, where $l'(\check{z}_i)$ is computed from an estimate of the marginal density $f$. Plugging this approximation into the likelihood factor in Proposition~\ref{prop:risk} yields an operational criterion.

\begin{proposition}\label{prop:approxRisk}
Given observations $\check{\pmb{z}}=(\check{z}_1,\ldots,\check{z}_d)^\top$, define $l'(\check{z}_i)=\frac{d}{d\check{z}_i}\log f(\check{z}_i)$ where $f$ is the marginal density of $\check{z}_i$. The risk of $\hat{\pmb{\theta}}_{\text{GEST}}$ can be approximated by
\[
\hat{R}(\hat{\pmb{\theta}}_{\text{GEST}}, \pmb{\check{z}})
=
d+\int_{\mathbb{R}^d}
\hat{H}\!\left(\pmb{z}, t, c, \frac{\partial t}{\partial |\pmb{z}|}, \frac{\partial c}{\partial |\pmb{z}|}\right)\,d\pmb{z},
\]
where
\begin{eqnarray*}
\hat{H}\!\left(\pmb{z}, t, c, \frac{\partial t}{\partial |\pmb{z}|}, \frac{\partial c}{\partial |\pmb{z}|}\right)
&=&
-\sum_{i=1}^d \frac{1}{|z_i|^{c(|\pmb{z}|)}}\Bigl(
2|z_i|\frac{\partial t(|\pmb{z}|)}{\partial |z_i|}
-
2t(|\pmb{z}|)|z_i|\ln|z_i|\frac{\partial c(|\pmb{z}|)}{\partial |z_i|}
\\
&&
+
2t(|\pmb{z}|)(1-c(|\pmb{z}|))
-
t(|\pmb{z}|)^2|z_i|^{2-c(|\pmb{z}|)}
\Bigr)
\prod_{j=1}^d \phi\!\left(z_j-\check{z}_j-l'(\check{z}_j)\right).
\end{eqnarray*}
\end{proposition}

\begin{remark}
We refer to $\hat{R}(\hat{\pmb{\theta}}_{\text{GEST}}, \pmb{\check{z}})$ in Proposition~\ref{prop:approxRisk} as the \emph{approximate risk}. It is a data-dependent criterion because it is constructed from $\check{\pmb{z}}$ and from an estimator of the marginal density $f$. It approximates the oracle risk by replacing the unknown $\pmb{\theta}$ inside the likelihood factor with a nonparametric plug-in approximation. Tweedie's formula is an empirical Bayes posterior mean and it is closely related to shrinkage rules such as James--Stein under Gaussian priors \cite{efron2011tweedie}. In our development, Tweedie's formula is used as an instrument for approximating the infeasible oracle risk of $\hat{\pmb{\theta}}_{\text{GEST}}$. Our goal is to minimize frequentist risk and its approximation over the GEST class rather than to minimize a Bayes risk over a prior family.
\end{remark}

\section{Approximate Risk Minimization and the NOMAD Estimator}\label{sec:nomad}

This section defines NOMAD in the canonical model as the estimator induced by minimizing the approximate risk over the GEST class. We will characterize minimizers of the approximate risk criterion through optimality conditions. 

\subsection{Approximate risk minimization over GEST}\label{sec:arm_over_gest}

\paragraph*{Objective and feasible set}
We consider the untruncated GEST family 
and an admissible class $\mathcal{F}$ of functionals $(t,c)$ such that the partial derivatives 
exist almost everywhere and the integrability conditions in Proposition~\ref{prop:risk} hold. 
For $(t,c)\in\mathcal{F}$, let $R(t,c;\pmb{\theta})$ denote the oracle risk in Proposition~\ref{prop:risk} and let $\hat{R}(t,c;\check{\pmb{z}})$ denote the approximate risk in Proposition~\ref{prop:approxRisk}.

\paragraph*{Approximate risk minimization and NOMAD}
We define $(\tilde t,\tilde c)$ as any minimizer of the approximate risk criterion,
\begin{equation}
(\tilde t,\tilde c)\in\arg\min_{(t,c)\in\mathcal{F}} \ \hat{R}(t,c;\check{\pmb{z}}),
\label{eq:arm_problem_hatR}
\end{equation}
and we define the resulting estimator by
\begin{equation}
\hat{\pmb{\theta}}_{\text{NOMAD}}
=
\check{\pmb{z}}+\pmb{h}(\check{\pmb{z}};\tilde t,\tilde c), \text{ where }
\pmb{h}(\check{\pmb{z}};\tilde t,\tilde c)
=-\tilde t(|\check{\pmb{z}}|)\Bigl(\sgn(\check{\pmb{z}})\circ|\check{\pmb{z}}|^{\,1-\tilde c(|\check{\pmb{z}}|)}\Bigr).
\label{eq:nomad_def_hatR}
\end{equation}
We refer to \eqref{eq:nomad_def_hatR} as the NOMAD estimator, which is short for \emph{NOnparametric Minimizer of an Approximate Data-driven Risk criterion}.

\subsection{Optimality conditions}\label{sec:arm_optimality}

We next record optimality conditions that characterize minimizers of the oracle risk and of the approximate risk. These conditions are derived by calculus of variations applied to the risk functional in Proposition~\ref{prop:risk} and to the approximate risk functional in Proposition~\ref{prop:approxRisk}. We defer the proof to our Supplement \cite{jiang2026nomadS}.

\subsubsection*{Oracle optimality conditions}
\begin{theorem}\label{thm:optRisk}
The functionals $t(|\pmb{z}|)$ and $c(|\pmb{z}|)$ minimize $R(\hat{\pmb{\theta}}_{\text{GEST}},\pmb{\theta})$ if and only if they satisfy
\begin{equation}
t(|\pmb{z}|)
=
\frac{\sum_{i=1}^d \sgn(z_i)|z_i|^{1-c(|\pmb{z}|)}(z_i-\theta_i)}
{\sum_{i=1}^d |z_i|^{2-2c(|\pmb{z}|)}},
\label{eq:riskCond1}
\end{equation}
and
\begin{equation}
\frac{\sum_{i=1}^d \sgn(z_i)|z_i|^{1-c(|\pmb{z}|)}(z_i-\theta_i)}
{\sum_{i=1}^d |z_i|^{2-2c(|\pmb{z}|)}}
=
\frac{\sum_{i=1}^d \sgn(z_i)|z_i|^{1-c(|\pmb{z}|)}\ln|z_i|\,(z_i-\theta_i)}
{\sum_{i=1}^d |z_i|^{2-2c(|\pmb{z}|)}\ln|z_i|}.
\label{eq:riskCond2}
\end{equation}
\end{theorem}

\subsubsection*{Approximate-risk optimality conditions}
\begin{theorem}\label{thm:optApproxRisk}
Let $\check{\pmb{z}}$ denote the observations and let $l'(\check{z}_i)=\frac{d}{d\check{z}_i}\log f(\check{z}_i)$. The functionals $t(|\pmb{z}|)$ and $c(|\pmb{z}|)$ minimize $\hat{R}(\hat{\pmb{\theta}}_{\text{GEST}},\pmb{\theta})$ if and only if they satisfy
\begin{equation}
t(|\check{\pmb{z}}|)
=
-\frac{\sum_{i=1}^d \sgn(\check{z}_i)|\check{z}_i|^{1-c(|\check{\pmb{z}}|)}\,l'(\check{z}_i)}
{\sum_{i=1}^d |\check{z}_i|^{2-2c(|\check{\pmb{z}}|)}},
\label{eq:approxRiskCond1}
\end{equation}
and
\begin{equation}
\frac{\sum_{i=1}^d \sgn(\check{z}_i)|\check{z}_i|^{1-c(|\check{\pmb{z}}|)}\,l'(\check{z}_i)}
{\sum_{i=1}^d |\check{z}_i|^{2-2c(|\check{\pmb{z}}|)}}
=
\frac{\sum_{i=1}^d \sgn(\check{z}_i)|\check{z}_i|^{1-c(|\check{\pmb{z}}|)}\ln|\check{z}_i|\,l'(\check{z}_i)}
{\sum_{i=1}^d |\check{z}_i|^{2-2c(|\check{\pmb{z}}|)}\ln|\check{z}_i|}.
\label{eq:approxRiskCond2}
\end{equation}
\end{theorem}

\paragraph*{Solving for $c$ and $t$}
Equation \eqref{eq:approxRiskCond2} yields an implicit equation for $c(|\check{\pmb{z}}|)$. Define
\begin{equation}
Q\!\left(c(|\check{\pmb{z}}|)\right):=AD-BC,
\label{eq:Qfunc}
\end{equation}
where
\begin{eqnarray*}
A&=&\sum_{i=1}^d \sgn(\check{z}_i)|\check{z}_i|^{1-c(|\check{\pmb{z}}|)}\,l'(\check{z}_i),\qquad
B=\sum_{i=1}^d |\check{z}_i|^{2-2c(|\check{\pmb{z}}|)},\\
C&=&\sum_{i=1}^d \sgn(\check{z}_i)|\check{z}_i|^{1-c(|\check{\pmb{z}}|)}\ln|\check{z}_i|\,l'(\check{z}_i),\qquad
D=\sum_{i=1}^d |\check{z}_i|^{2-2c(|\check{\pmb{z}}|)}\ln|\check{z}_i|.
\end{eqnarray*}
We solve $Q(c)=0$ by Newton iteration. Given $c(|\check{\pmb{z}}|)$, we compute $t(|\check{\pmb{z}}|)$ from \eqref{eq:approxRiskCond1}.

\section{Connections with Classical Estimators}\label{sec:connections}

This section clarifies how classical estimators arise from the approximate-risk optimality conditions in Theorem~\ref{thm:optApproxRisk}. The central observation is that these conditions depend on the marginal score function $l'(z)$, 
where $f$ denotes the coordinatewise marginal density of the empirical observations. Different shapes of $f$ induce different shrinkage-thresholding profiles through the pair $(\tilde t,\tilde c)$. As a result, familiar procedures such as maximum likelihood estimation, James--Stein shrinkage, ridge-type shrinkage, and lasso-type soft-thresholding appear as exact or limiting solutions under specific score structures \cite{efron1973stein,park2008bayesian}. We defer the proof of propositions to Supplement \cite{jiang2026nomadS} 

\subsection{Reduction to maximum likelihood estimation}\label{sec:conn_mle}

\paragraph*{Vanishing score implies no shrinkage}
If the marginal score is identically zero, then the approximate-risk criterion provides no incentive to shrink the observations. In this case, the optimality condition \eqref{eq:approxRiskCond1} immediately gives $\tilde t(|\check{\pmb z}|)=0$, so the resulting rule coincides with the maximum likelihood estimator.

\begin{proposition}[Reduction to MLE]\label{prop:rel_mle}
If the marginal score satisfies $l'(\check z_i)=0$ for all $i=1,\ldots,d$, then the approximate-risk minimizer satisfies $\tilde t(|\check{\pmb z}|)=0$, and the resulting NOMAD estimator equals the maximum likelihood estimator,
\[
\hat{\pmb{\theta}}_{\mathrm{NOMAD}}=\check{\pmb z}.
\]
\end{proposition}


\subsection{Reduction to James--Stein and ridge-type shrinkage}\label{sec:conn_js_ridge}

\paragraph*{Gaussian marginals yield linear score}
Suppose the marginal density is centered Gaussian with variance $\tau^2$. Then
\begin{equation}
l'(z)=-\frac{z}{\tau^2}.
\label{eq:linear_score}
\end{equation}
A linear score favors shrinkage that is proportional to the observed magnitude and does not favor thresholding. Substituting \eqref{eq:linear_score} into the optimality conditions shows that the minimizing profile satisfies $c(\cdot)=0$, so the resulting rule is purely shrinkage based.

\begin{proposition}[Gaussian marginal and James--Stein reduction]\label{prop:rel_js}
Suppose the empirical marginal density $f$ is $N(0,\tau^2)$, so that $l'(z)=-z/\tau^2$. Then the approximate-risk minimizer yields a shrinkage-only rule with
\[
c(|\check{\pmb z}|)=0,
\qquad
\hat{\pmb{\theta}}_{\mathrm{NOMAD}}
=
\left(1-\tau^{-2}\right)\check{\pmb z}.
\]
Under the usual empirical-Bayes substitution of the unknown shrinkage level by its James--Stein estimate, this profile reduces to the James--Stein rule
\[
\hat{\pmb{\theta}}_{\mathrm{JS}}
=
\left(1-\frac{d-2}{\|\check{\pmb z}\|_2^2}\right)\check{\pmb z}.
\]
Equivalently, the induced rule belongs to the ridge-type shrinkage family.
\end{proposition}


\subsection{Reduction to lasso and thresholding-type estimators}\label{sec:conn_lasso}

\paragraph*{Laplace-like tails yield nearly constant score magnitude}
Soft-thresholding behavior arises when the marginal density is symmetric, sharply peaked near zero, and has tails whose score is approximately constant in magnitude. More precisely, suppose there exists $\kappa>0$ and $u>0$ such that
\begin{equation}
l'(z)\approx -\kappa\,\sgn(z)
\qquad\text{for all } |z|\ge u.
\label{eq:laplace_like_score}
\end{equation}
This condition holds exactly for a centered Laplace density. In the present setting, however, the empirical marginal law is typically a Gaussian location mixture, so \eqref{eq:laplace_like_score} should be interpreted as an effective tail-shape condition rather than a literal model assumption \cite{johnstone2005empirical}.

\begin{proposition}[Laplace marginal and lasso reduction]\label{prop:rel_lasso}
If the marginal distribution $f(Z)$ is a centered Laplace distribution with score
\[
l'(z)=-\kappa\,\sgn(z),
\]
then $\hat R(\hat{\pmb{\theta}}_{\mathrm{GEST}},\pmb{\theta})$ is minimized if and only if
\[
t(|\pmb Z|)=\kappa
\qquad\text{and}\qquad
c(|\pmb Z|)=1.
\]
In this case, the realized estimator is the lasso soft-thresholding rule with threshold $\lambda=\kappa$, namely
\[
\hat{\pmb{\theta}}_{\mathrm{NOMAD}}
=
\hat{\pmb{\theta}}_{\mathrm{lasso}}
=
\left(1-\frac{\lambda}{|\pmb Z|}\right)_+\circ \pmb Z.
\]
\end{proposition}

When \eqref{eq:laplace_like_score} holds over the range of observed coordinates, the numerator in \eqref{eq:approxRiskCond1} behaves like a signed weighted $\ell_1$ term. The optimality conditions then favor choices with $c(|\check{\pmb z}|)$ close to one. In this regime, the induced shrinkage multiplier takes the soft-thresholding form, and after positive-part truncation the resulting estimator produces exact zeros. Thus the approximate-risk minimizer reduces to the lasso rule \cite{park2008bayesian}.

\paragraph*{Intermediate sparse regimes}
If the empirical marginal density is sharply peaked near zero and heavy tailed, but the score is not exactly constant in magnitude, then the approximate-risk minimizer typically selects
\[
0<c(|\check{\pmb z}|)<1.
\]
This produces mixed shrinkage-thresholding behavior, with stronger attenuation for small coordinates and milder shrinkage for large ones. This intermediate regime explains how NOMAD interpolates between James--Stein-type shrinkage in dense settings and lasso-type thresholding in sparse settings \cite{johnstone2005empirical,abramovich2006adapting}.

\begin{proposition}\label{prop:range}
If the marginal distribution $f(Z)$ is an even function and monotonically decreasing in $|Z|$, then $\hat R(\hat{\pmb{\theta}}_{\mathrm{GEST}},\pmb{\theta})$ can be minimized only if
\[
0\le c(|\pmb z|)\le 1.
\]
\end{proposition}



\section{Approximate Risk Minimization Theory in the Canonical Model}\label{sec:canonical_theory}

This section develops a consistency theory for NOMAD in the canonical normal mean model under a structured admissible class. The main conclusion is that, under regularity conditions on the signal sequence, the score estimator, and the sieve class, the approximate risk uniformly tracks the oracle risk over the sieve and the resulting NOMAD estimator is oracle-risk consistent on the scale of the dimension. Under an additional identifiability condition, the selected shrinkage-thresholding functionals are also consistent.

\subsection{Primitive assumptions}\label{sec:canonical_assumptions}

We organize the primitive assumptions into three groups. The first group specifies the canonical model and the structured admissible class. The second group gives a concrete sieve construction. The third group imposes signal-distribution and score-estimation conditions.

\paragraph*{(A) Canonical model and structured admissible function class}
Assume
\[
\check z_i=\theta_i+\varepsilon_i,
\qquad
\varepsilon_i\overset{\mathrm{iid}}{\sim}N(0,1),
\qquad
i=1,\ldots,d,
\]
with deterministic \(\pmb\theta=\pmb\theta^{(d)}\). Let
\[
G_d=\frac1d\sum_{i=1}^d \delta_{\theta_i}
\]
denote the empirical signal distribution.

Let \(q_t,q_c\ge 1\) be fixed integers. For \(m=1,\ldots,q_t\) and \(\ell=1,\ldots,q_c\), let \(\psi_{tm}\) and \(\psi_{c\ell}\) be \(C^2\) functions on \(\mathbb R\) such that, for some constants \(C_\psi>0\) and \(r_\psi\ge 0\),
\[
|\psi(x)|+|\psi'(x)|+|\psi''(x)|
\le
C_\psi(1+|x|^{r_\psi})
\qquad\text{for all }x\in\mathbb R.
\]
Define the summary maps
\[
s_t(\pmb z)
=
\left(
\frac1d\sum_{i=1}^d \psi_{t1}(z_i),\ldots,\frac1d\sum_{i=1}^d \psi_{tq_t}(z_i)
\right)\in\mathbb R^{q_t},
\]
\[
s_c(\pmb z)
=
\left(
\frac1d\sum_{i=1}^d \psi_{c1}(z_i),\ldots,\frac1d\sum_{i=1}^d \psi_{cq_c}(z_i)
\right)\in\mathbb R^{q_c}.
\]
Let \(\mathcal S_t\subset\mathbb R^{q_t}\) and \(\mathcal S_c\subset\mathbb R^{q_c}\) be compact sets containing the ranges of these summaries with overwhelming probability.

We consider the admissible class $\mathcal F$ for functionals of the form
\[
t(\pmb z)=\tau(s_t(\pmb z)),
\qquad
c(\pmb z)=\gamma(s_c(\pmb z)),
\]
where
\[
\tau:\mathcal S_t\to[0,t_{\max}],
\qquad
\gamma:\mathcal S_c\to[0,c_{\max}]
\]
with \(0\le c_{\max}<1\), and \(\tau,\gamma\in C^2\) with uniformly bounded first and second derivatives:
\[
\sup_{u\in\mathcal S_t}\|\nabla\tau(u)\|_2
+
\sup_{u\in\mathcal S_t}\|D^2\tau(u)\|_{\mathrm{op}}
\le L_t,
\]
\[
\sup_{v\in\mathcal S_c}\|\nabla\gamma(v)\|_2
+
\sup_{v\in\mathcal S_c}\|D^2\gamma(v)\|_{\mathrm{op}}
\le L_c.
\]

For \((\tau,\gamma)\in\mathcal F\), let \((t,c)\) denote the induced functionals, let \(\hat{\pmb\theta}_{t,c}\) denote the corresponding GEST estimator, and define the oracle quadratic risk
\[
R(t,c;\pmb\theta)
=
E_{\pmb\theta}\!\left[\|\hat{\pmb\theta}_{t,c}-\pmb\theta\|_2^2\right].
\]

\paragraph*{(B) Concrete sieve construction}
Let \(\mathcal T_d\) and \(\mathcal C_d\) be tensor-product spline sieves of order at least \(3\) on \(\mathcal S_t\) and \(\mathcal S_c\), respectively. Write
\[
K_{t,d}=\dim(\mathcal T_d),
\qquad
K_{c,d}=\dim(\mathcal C_d),
\qquad
m_d=K_{t,d}+K_{c,d}.
\]
Assume
\[
K_{t,d}\to\infty,
\qquad
K_{c,d}\to\infty,
\qquad
m_d\to\infty,
\qquad
\frac{m_d}{d}\to 0.
\]
The spline coefficient constraints are chosen so that every element of \(\mathcal T_d\) remains in \([0,t_{\max}]\) and has first and second derivatives bounded by \(L_t\), and every element of \(\mathcal C_d\) remains in \([0,c_{\max}]\) and has first and second derivatives bounded by \(L_c\). Define
\[
\mathcal F_d=\mathcal T_d\times\mathcal C_d\subset\mathcal F.
\]

We use the \(C^1\)-type metric
\[
d_{\mathcal F}^{(1)}\bigl((\tau_1,\gamma_1),(\tau_2,\gamma_2)\bigr)
=
\|\tau_1-\tau_2\|_{\infty,\mathcal S_t}
+
\|\nabla\tau_1-\nabla\tau_2\|_{\infty,\mathcal S_t}
+
\|\gamma_1-\gamma_2\|_{\infty,\mathcal S_c}
+
\|\nabla\gamma_1-\nabla\gamma_2\|_{\infty,\mathcal S_c}.
\]
Assume the spline sieves are dense in \(\mathcal F\) under \(d_{\mathcal F}^{(1)}\), that is, for every \((\tau,\gamma)\in\mathcal F\) there exists \((\tau_d,\gamma_d)\in\mathcal F_d\) such that
\[
d_{\mathcal F}^{(1)}\bigl((\tau_d,\gamma_d),(\tau,\gamma)\bigr)\to 0.
\]

\paragraph*{(C) Signal tails, score estimation, and oracle Tweedie accuracy}
Assume \(G_d\rightarrow G\) for some probability measure \(G\) on \(\mathbb R\), and assume a uniform sub-exponential tail condition:
\begin{equation}
\sup_d \frac1d\sum_{i=1}^d e^{\lambda_0 |\theta_i|}<\infty
\label{eq:subexp_signal}
\end{equation}
for some \(\lambda_0>0\).

Let
$
f_d(z)=\int \phi(z-\theta)\,dG_d(\theta)$,
$f(z)=\int \phi(z-\theta)\,dG(\theta)$,
and define
\[
l_d'(z)=\frac{d}{dz}\log f_d(z),
\qquad
l'(z)=\frac{d}{dz}\log f(z).
\]
Let \(\widehat l_d'\) be a score estimator obtained from a regular density estimator, for example KDE with derivative or a spline estimator of the log-density. Assume there exists a deterministic sequence \(a_d\to\infty\) such that
\[
\sup_{|z|\le a_d}|\widehat l_d'(z)-l'(z)|\overset{p}\to 0
\text{ and }
\sup_{|z|\le a_d}|l_d'(z)-l'(z)|\to 0.
\]
Assume also that there exists a measurable envelope \(M_d(z)\) such that
\[
|\widehat l_d'(z)|+|l_d'(z)|\le M_d(z)
\qquad\text{for all }z,
\]
and for some fixed \(q\ge 2\),
\[
M_d(z)^2\le C_M(1+|z|^q)
\]
uniformly in \(d\). Finally, assume the oracle Tweedie posterior mean is accurate on average:
\begin{equation}
\frac1d\sum_{i=1}^d
E\!\left[
\bigl(\theta_i-\check z_i-l_d'(\check z_i)\bigr)^2
\right]\to 0.
\label{eq:oracle_tweedie_accuracy_final}
\end{equation}

\subsection{Oracle-risk consistency of NOMAD}\label{sec:oracle_consistency}


Define the oracle value
\[
R^\star(\pmb\theta)=\inf_{(\tau,\gamma)\in\mathcal F}R(t,c;\pmb\theta)
\]
and the oracle minimizer set
\[
\mathcal A^\star(\pmb\theta)=\arg\min_{(\tau,\gamma)\in\mathcal F}R(t,c;\pmb\theta),
\]
where \((t,c)\) is induced by \((\tau,\gamma)\). By compactness of \(\mathcal F\) under \(d_{\mathcal F}^{(1)}\) and continuity of risk (See Supplement \cite{jiang2026nomadS}) 
, \(\mathcal A^\star(\pmb\theta)\) is nonempty. Fix a measurable oracle selection \((\tau^\star,\gamma^\star)\in\mathcal A^\star(\pmb\theta)\).

For each \(d\), define the sieve NOMAD pair by
\begin{equation}
(\tilde\tau_d,\tilde\gamma_d)\in\arg\min_{(\tau,\gamma)\in\mathcal F_d}\hat R(t,c;\check{\pmb z}),
\label{eq:sieve_nomad_pair_final}
\end{equation}
with a deterministic tie-breaking rule if needed, and define
\[
\hat{\pmb\theta}_{\mathrm{NOMAD}}=\hat{\pmb\theta}_{\tilde t_d,\tilde c_d},
\]
where \((\tilde t_d,\tilde c_d)\) is induced by \((\tilde\tau_d,\tilde\gamma_d)\).

The next theorem states the main consistency result for the canonical model.

\begin{theorem}[Oracle-risk consistency relative to \(\mathcal F\)]
\label{thm:oracle_consistency}
Assume (A), (B), (C), then
\[
R(\hat{\pmb\theta}_{\mathrm{NOMAD}},\pmb\theta)-R^\star(\pmb\theta)
=
o_{\mathbb P}(d).
\]
\end{theorem}


We finally state functional consistency under an identifiability condition.

\paragraph*{Identifiability condition}
Assume the oracle minimizer \((\tau^\star,\gamma^\star)\) is unique and there exists a strictly increasing function \(\varphi:[0,\infty)\to[0,\infty)\) with \(\varphi(0)=0\) such that
\[
\frac1d\bigl(R(t,c;\pmb\theta)-R^\star(\pmb\theta)\bigr)
\ge
\varphi\!\left(d_{\mathcal F}^{(1)}\bigl((\tau,\gamma),(\tau^\star,\gamma^\star)\bigr)\right)
\]
for all \((\tau,\gamma)\in\mathcal F\).

\begin{theorem}[Functional consistency]
\label{thm:functional_consistency}
Assume the conditions of Theorem~\ref{thm:oracle_consistency} and the identifiability condition above. Then
\[
d_{\mathcal F}^{(1)}\bigl((\tilde\tau_d,\tilde\gamma_d),(\tau^\star,\gamma^\star)\bigr)\overset{p}\to 0.
\]
\end{theorem}

We defer the proof of both theorems to the Supplement \cite{jiang2026nomadS}. 

\subsection{Approximate-risk inequality}\label{sec:approx_risk_inequality}

The minimization defining NOMAD yields a deterministic ordering on the approximate-risk scale.

\begin{theorem}[Approximate-risk inequality]
\label{thm:approx_risk_inequality}
Let \((\tilde\tau_d,\tilde\gamma_d)\) satisfy \eqref{eq:sieve_nomad_pair_final}. Then for any \((\tau,\gamma)\in\mathcal F_d\),
\[
\hat R(\tilde t_d,\tilde c_d;\check{\pmb z})
\le
\hat R(t,c;\check{\pmb z}),
\]
where \((\tilde t_d,\tilde c_d)\) and \((t,c)\) are induced by \((\tilde\tau_d,\tilde\gamma_d)\) and \((\tau,\gamma)\), respectively.
\end{theorem}

\section{Application to Wavelet Analysis}\label{sec:wavelets}

This section applies the canonical framework to wavelet denoising. The wavelet transform converts a signal-plus-noise problem to coefficient-wise normal mean estimation under standard assumptions, so the GEST class and the approximate risk minimization principle can be applied directly \cite{donoho1994ideal,donoho1995adapting}. The application provides an interpretable illustration of how NOMAD adapts between shrinkage and thresholding in a classical nonparametric setting.

\subsection{Wavelet denoising as a normal mean problem}\label{sec:wavelets_model}

\paragraph*{Signal-plus-noise model}
Let $\pmb{y}=(y_1,\ldots,y_n)^\top$ denote observations on an equispaced grid and consider the Gaussian sequence model
\begin{equation}
y_i = s_i + \sigma \varepsilon_i,
\qquad
i=1,\ldots,n,
\qquad
\varepsilon_i \overset{\mathrm{iid}}{\sim} N(0,1),
\label{eq:wavelet_signal_noise}
\end{equation}
where $\pmb{s}=(s_1,\ldots,s_n)^\top$ is the unknown signal and $\sigma>0$ is the noise level. We treat $\sigma$ as known for exposition, and in practice it can be estimated from the finest-scale coefficients by a robust median-based estimator \cite{donoho1995adapting}.

\paragraph*{Wavelet transform and coefficient representation}
Let $\pmb{W}$ denote an orthonormal discrete wavelet transform matrix, so $\pmb{W}^\top \pmb{W}=\pmb{I}_n$. Define wavelet coefficients $\pmb{Z}=\pmb{W}\pmb{y}$ and $\pmb{\theta}=\pmb{W}\pmb{s}$. Multiplying \eqref{eq:wavelet_signal_noise} by $\pmb{W}$ yields
\begin{equation}
\pmb{Z} = \pmb{\theta} + \sigma \pmb{\xi},
\qquad
\pmb{\xi}\sim N(\pmb{0},\pmb{I}_n),
\label{eq:wavelet_normal_means}
\end{equation}
which is a multivariate normal mean model. After rescaling $\pmb{Z}/\sigma$, \eqref{eq:wavelet_normal_means} matches the canonical form \eqref{eq:canonical_model}. This reduction is exact for orthonormal transforms and independent Gaussian noise, and it is approximately valid when the transform is nearly orthogonal \cite{mallat1989theory,donoho1994ideal}.

\paragraph*{Role of thresholding in wavelet denoising}
In wavelet denoising, shrinkage and thresholding are applied to wavelet coefficients to suppress noise while retaining structure. Classical rules such as VisuShrink use a universal hard or soft threshold, and SureShrink uses Stein's unbiased risk estimate to select thresholds in a data-adaptive manner \cite{donoho1994ideal,donoho1995adapting}. These approaches motivate treating wavelet denoising as a selection problem over a family of shrinkage-thresholding rules, which matches the developed approximate risk minimization framework. 

\subsection{NOMAD in the wavelet domain}\label{sec:wavelets_nomad}

\paragraph*{Coefficient-domain NOMAD rule}
Let $\check{\pmb{z}}=\pmb{W}\pmb{y}/\hat{\sigma}$ denote standardized empirical coefficients where $\hat{\sigma}$ is a noise estimator. We apply the canonical NOMAD construction to $\check{\pmb{z}}$ and obtain a coefficient estimator
$\hat{\pmb{\theta}}^{\,w}_{\mathrm{NOMAD}}
=
\hat{\pmb{\theta}}_{\mathrm{NOMAD}}(\check{\pmb{z}}),$
where $\hat{\pmb{\theta}}_{\mathrm{NOMAD}}$ is defined in \eqref{eq:nomad_def_hatR} through approximate risk minimization over the GEST class. The denoised signal is then obtained by inverse transform,
\begin{equation}
\hat{\pmb{s}}_{\mathrm{NOMAD}}
=
\pmb{W}^\top \bigl(\hat{\sigma}\,\hat{\pmb{\theta}}^{\,w}_{\mathrm{NOMAD}}\bigr).
\label{eq:wavelet_inverse}
\end{equation}

\paragraph*{Level-dependent implementation}
Wavelet coefficients are naturally grouped by resolution levels. To respect this structure and to reduce sensitivity to heterogeneity across scales, we implement NOMAD in a level-dependent manner. Let $\mathcal{I}_j$ denote the index set of coefficients at level $j$, and let $\check{\pmb{z}}^{(j)}=\{\check{z}_i:i\in\mathcal{I}_j\}$. We apply approximate risk minimization to each level 
$
\hat{\pmb{\theta}}^{\,w,(j)}_{\mathrm{NOMAD}}
=
\hat{\pmb{\theta}}_{\mathrm{NOMAD}}(\check{\pmb{z}}^{(j)})
$
, and we concatenate the levelwise estimates to form $\hat{\pmb{\theta}}^{\,w}_{\mathrm{NOMAD}}$. This approach preserves the canonical normal means structure within each level and it allows the estimated shrinkage-thresholding profile to vary across levels \cite{donoho1995adapting,johnstone2005empirical}.

\paragraph*{Choice of decomposition level and boundary handling}
We use a standard dyadic decomposition with $J=\log_2(n)$ levels when $n$ is a power of two. For general $n$ we use standard padding or boundary wavelet schemes that yield a nearly orthonormal transform \cite{mallat1989theory}. 

\subsection{Comparative empirical results}\label{sec:wavelets_results}

\begin{figure}[h!]
    \centering
    \begin{subfigure}[b]{0.24\textwidth}
        \centering
        \includegraphics[width=\textwidth]{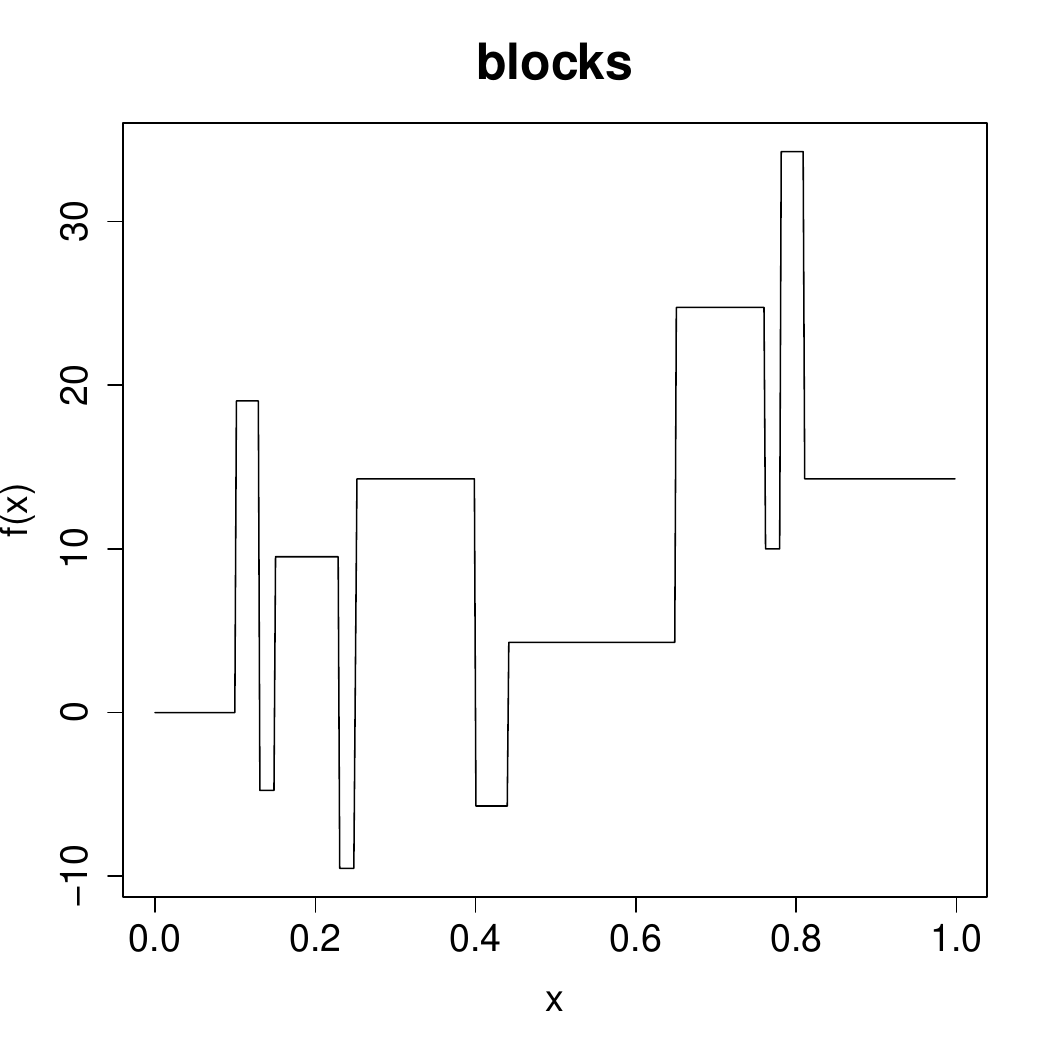}
    \end{subfigure}
    \begin{subfigure}[b]{0.24\textwidth}
        \centering
        \includegraphics[width=\textwidth]{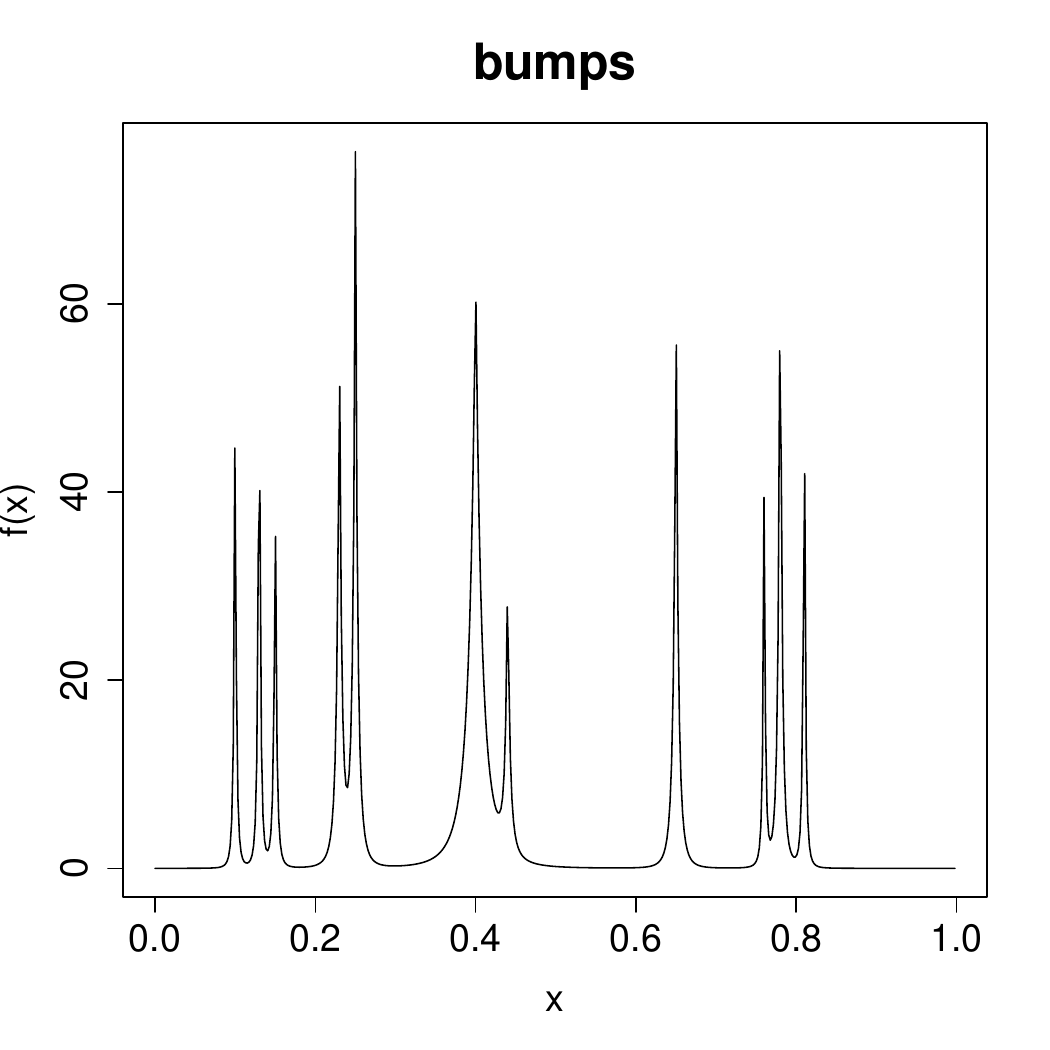}
    \end{subfigure}
    \begin{subfigure}[b]{0.24\textwidth}
        \centering
        \includegraphics[width=\textwidth]{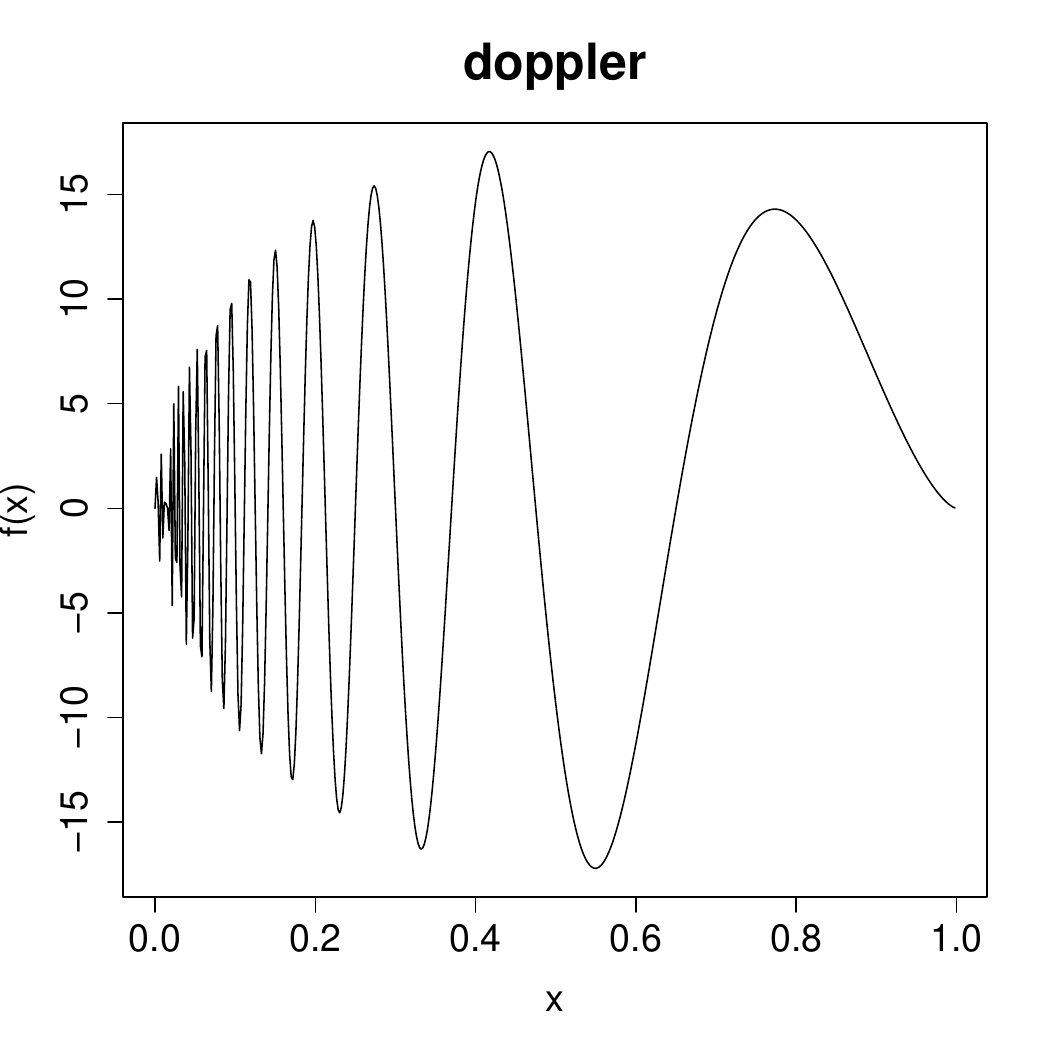}
    \end{subfigure}
    \begin{subfigure}[b]{0.24\textwidth}
        \centering
        \includegraphics[width=\textwidth]{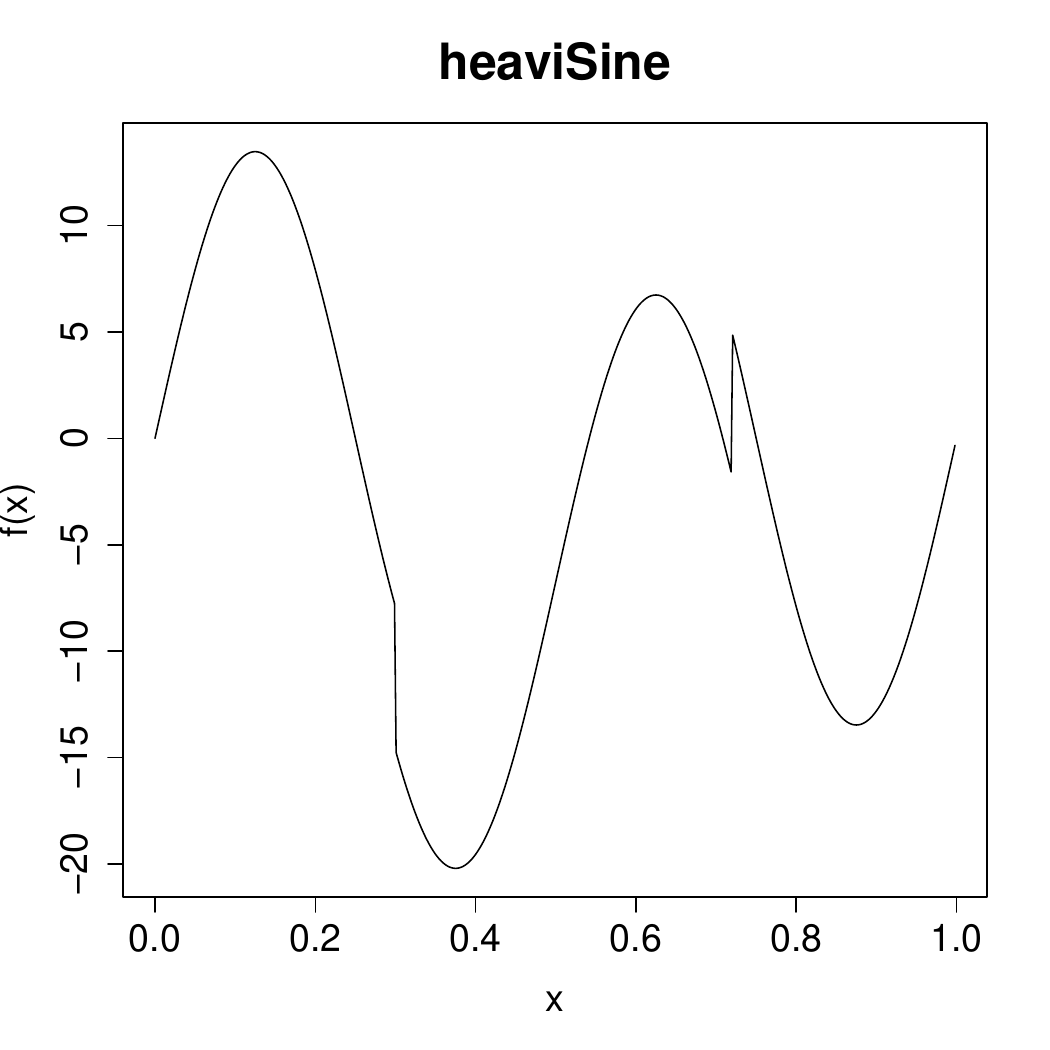}
    \end{subfigure}\\
    \caption{Representative signal shapes used in the wavelet denoising experiments.}
    \label{fig:wavelet_signal}
\end{figure}

We compare $\hat{\pmb{s}}_{\mathrm{NOMAD}}$ with standard wavelet denoising procedures, including universal soft-thresholding (VisuShrink) and levelwise SURE thresholding (SureShrink). VisuShrink applies a universal soft-threshold to the empirical wavelet coefficients, with the threshold scaled by an estimate of the noise level based on the median absolute deviation. SureShrink also uses soft-thresholding, but it selects the threshold by minimizing Stein's unbiased risk estimate level by level. These benchmarks represent two classical paradigms in wavelet denoising, namely fixed-threshold rules and risk-estimated threshold rules \cite{donoho1994ideal,donoho1995adapting}. In addition, we include the positive-part James--Stein estimator (JS+) in the comparison. For all methods, performance is evaluated by the mean squared error
$
\mathrm{MSE}(\hat{\pmb{s}})=\frac{1}{n}\|\hat{\pmb{s}}-\pmb{s}\|_2^2,
$
and results are reported across a range of signal classes and sample sizes in order to highlight regime-dependent behavior. The underlying signal shapes are displayed in Figure~\ref{fig:wavelet_signal}.

\begin{figure}[h!]
    \centering
    \begin{subfigure}[b]{0.4\textwidth}
        \centering
        \includegraphics[width=\textwidth]{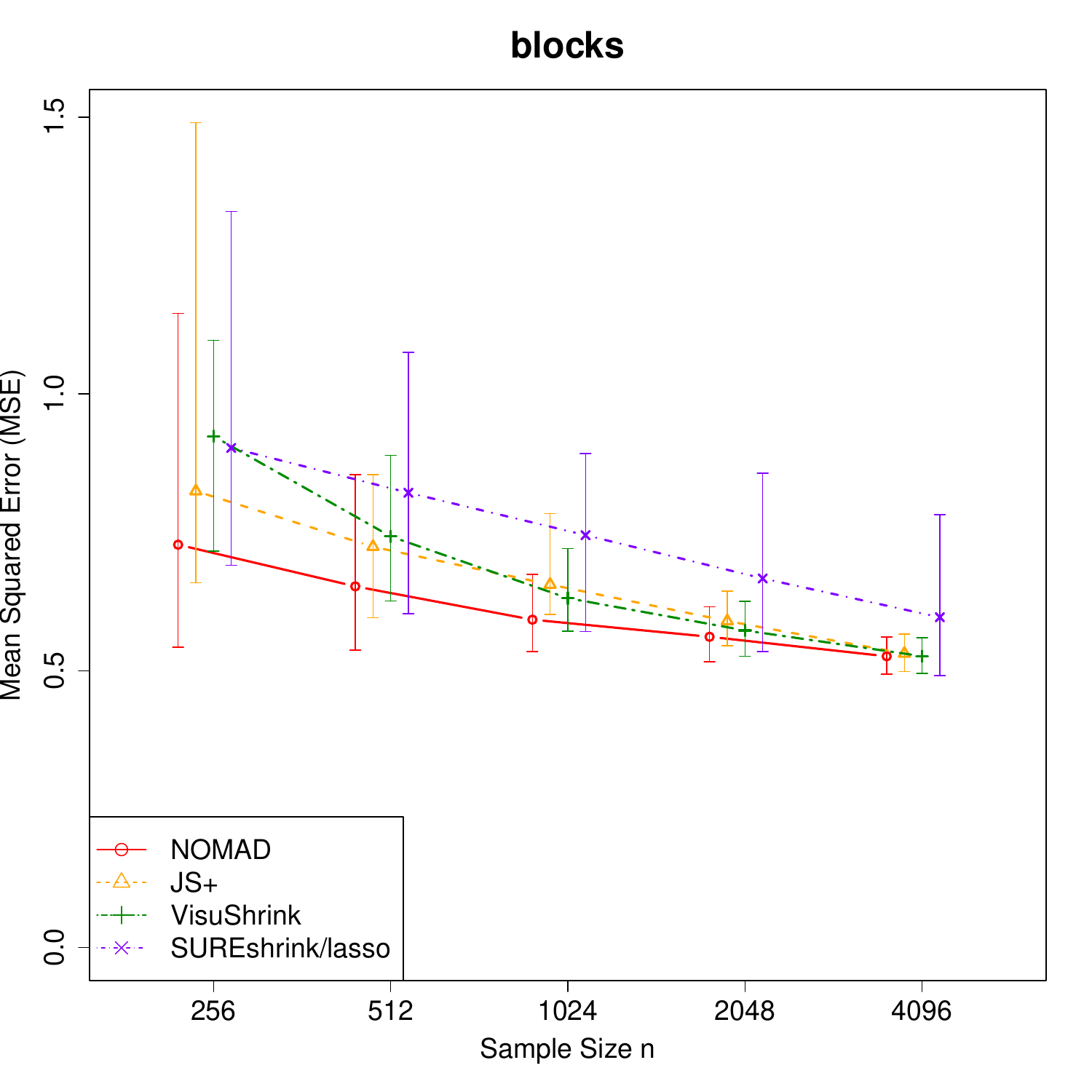}
    \end{subfigure}
    \begin{subfigure}[b]{0.4\textwidth}
        \centering
        \includegraphics[width=\textwidth]{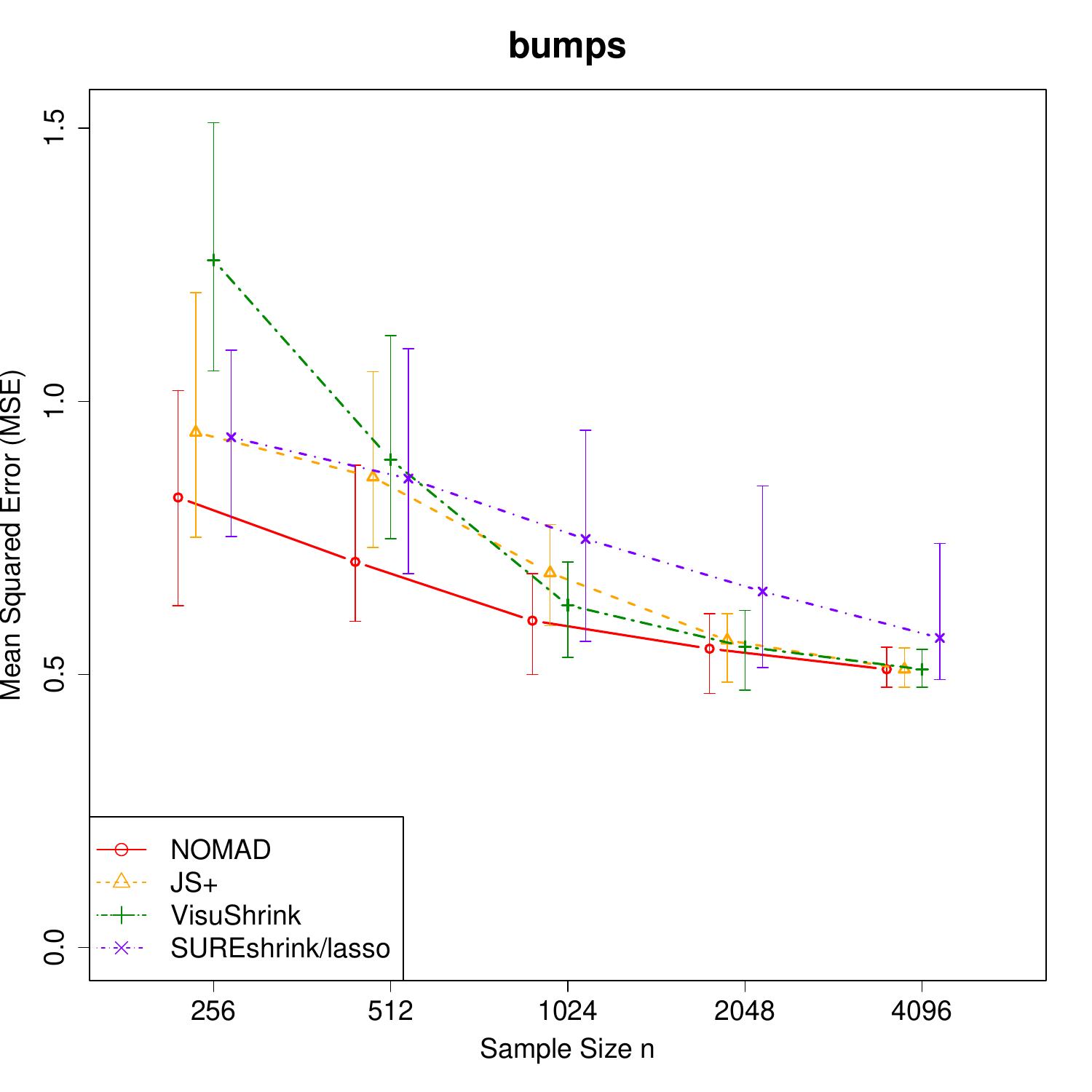}
    \end{subfigure}\\
    \begin{subfigure}[b]{0.4\textwidth}
        \centering
        \includegraphics[width=\textwidth]{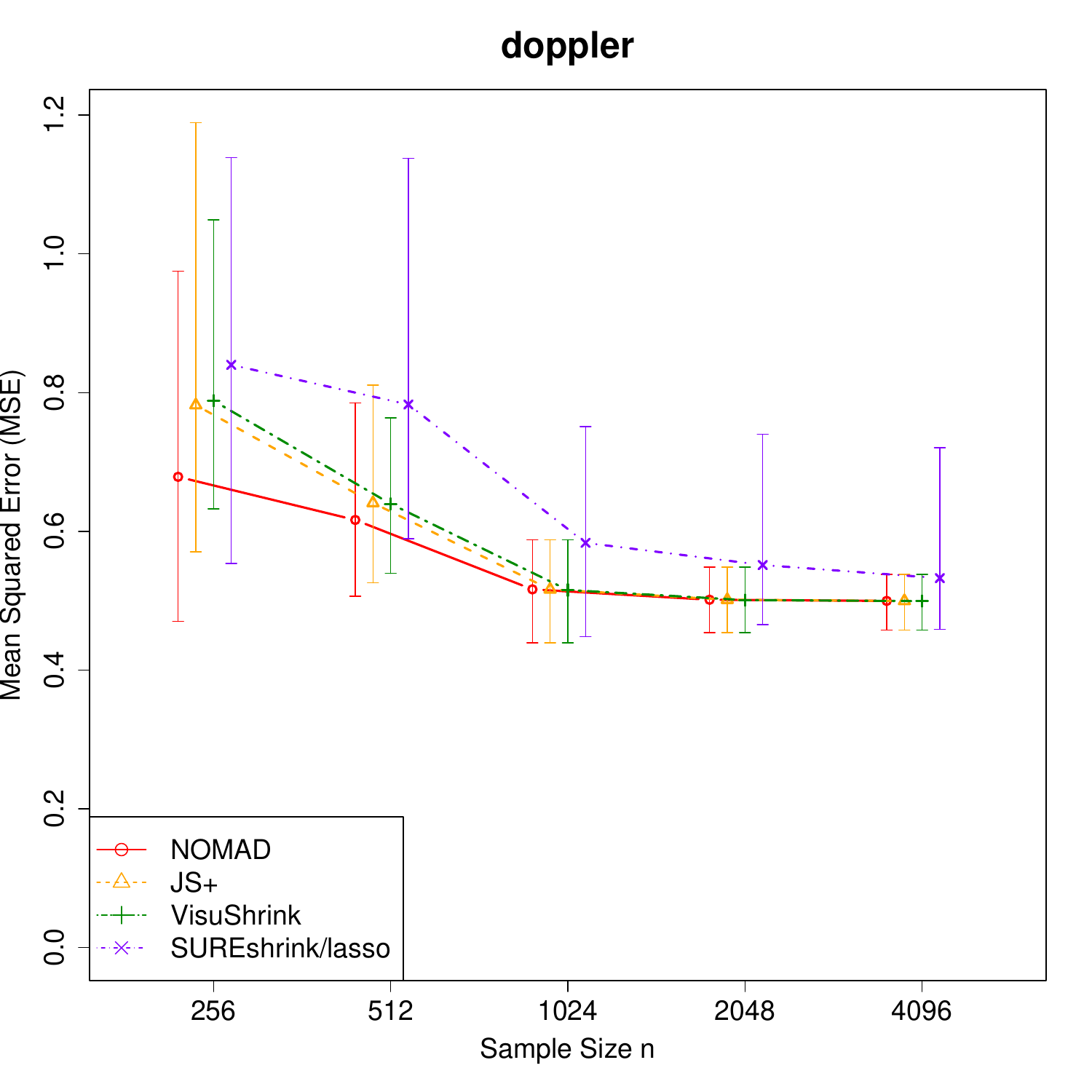}
    \end{subfigure}
    \begin{subfigure}[b]{0.4\textwidth}
        \centering
        \includegraphics[width=\textwidth]{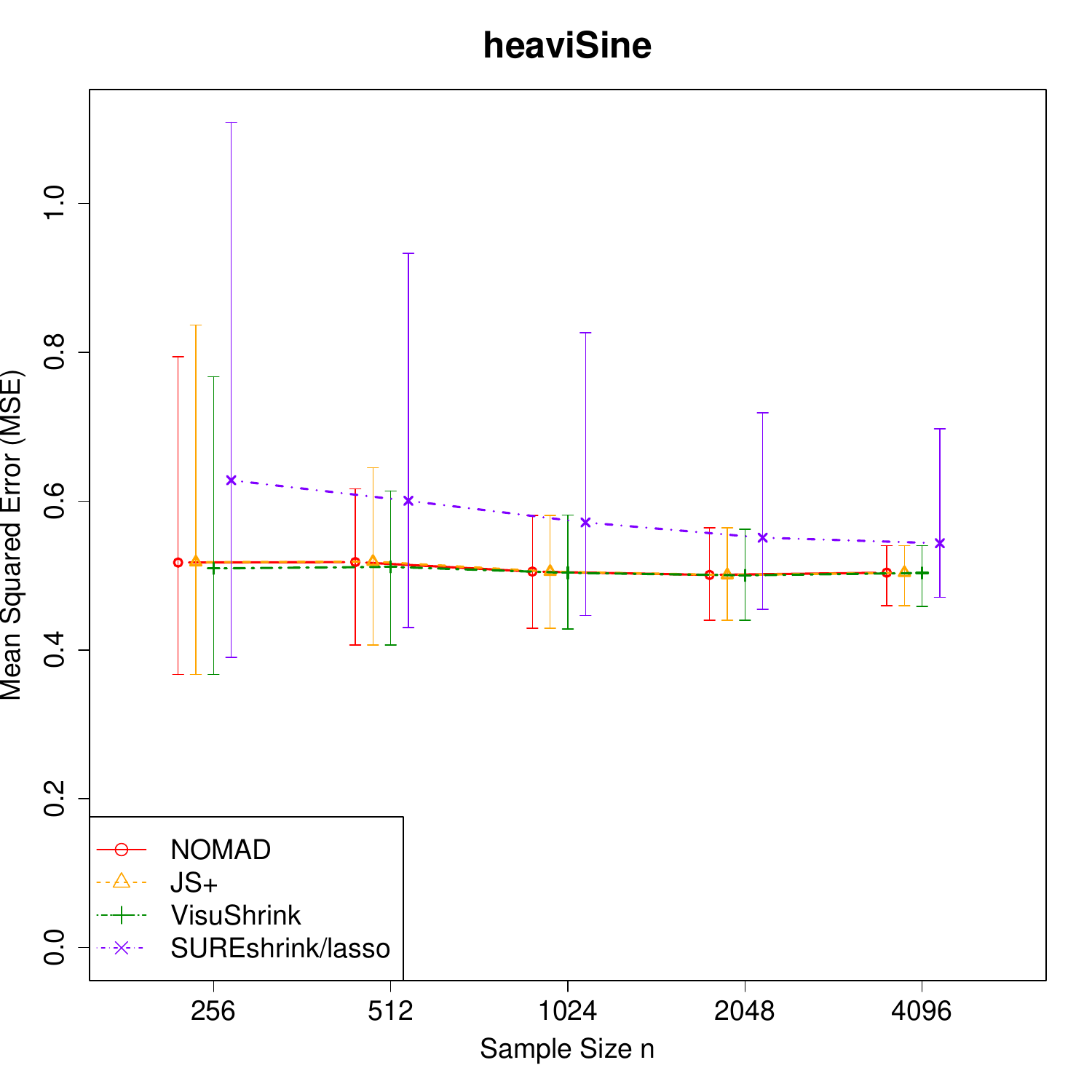}
    \end{subfigure}\\
    
    \caption{Mean squared error of competing wavelet denoising methods across signal classes and sample sizes.}
    \label{fig:wavelet_mse}
\end{figure}

Figure~\ref{fig:wavelet_mse} summarizes the denoising performance of the competing methods across the different signal classes and sample sizes. NOMAD attains the smallest mean squared error for the Doppler, Bumps, and Blocks signals. For HeaviSine, the performances of NOMAD, JS+, and VisuShrink are nearly indistinguishable. This behavior is consistent with the empirical marginal distribution of the wavelet coefficients in that example, which is close to a centered Gaussian distribution; see Supplement \cite{jiang2026nomadS}. In that regime, the theory developed in Section~\ref{sec:connections} predicts that the optimal exponent satisfies $c=0$, so the NOMAD rule reduces to a shrinkage-only procedure and behaves like the James--Stein estimator.

Wavelet coefficients exhibit substantial heterogeneity across levels and across signal classes. The NOMAD construction adapts to this heterogeneity through the score-driven approximate risk criterion, and it can interpolate between shrinkage-dominant behavior at coarse scales and thresholding-dominant behavior at fine scales. This behavior agrees with the classical intuition that coarse-scale coefficients are often relatively dense, whereas fine-scale coefficients are often sparse \cite{johnstone2005empirical}.

Wavelet denoising illustrates the framework in a setting that reduces to a canonical normal mean problem after an orthogonal transform \cite{donoho1994ideal,donoho1995adapting}. Many applications, however, involve correlated observations or more general dependence structures for which such a reduction is unavailable. The next section extends approximate risk minimization to multivariate normal mean estimation with correlated observations, and that extension in turn provides the foundation for the linear regression formulation developed later.

\section{Extension to Correlated Multivariate Normal Mean Estimation}\label{sec:correlated_extension}

This section extends the canonical framework to multivariate normal mean estimation with correlation. The extension has two purposes. First, it shows how shrinkage-thresholding can be defined when the covariance structure is nontrivial and the natural loss is the invariant quadratic risk induced by the precision matrix \cite{berger1976admissible}. Second, it prepares the regression formulation in Section~\ref{sec:regression}, where the Gram matrix plays the role of a precision matrix.

\subsection{Model setup and challenges under correlation}\label{sec:correlated_setup}

\paragraph*{General multivariate normal mean model}
Let
\[
\pmb{Z}\sim N(\pmb{\theta},\pmb{\Omega}^{-1}),
\]
where $\pmb{Z}=(Z_1,\ldots,Z_d)^\top$ is the observed vector, $\pmb{\theta}=(\theta_1,\ldots,\theta_d)^\top$ is the unknown mean vector, and $\pmb{\Omega}$ is a positive definite precision matrix. In this setting, a natural 
risk is the invariant quadratic risk
\[
R(\hat{\pmb{\theta}},\pmb{\theta})
=
E_{\pmb{\theta}}\!\left[(\hat{\pmb{\theta}}-\pmb{\theta})^\top \pmb{\Omega}(\hat{\pmb{\theta}}-\pmb{\theta})\right].
\]
When $\pmb{\Omega}=\pmb{I}_d$, this reduces to the canonical normal mean problem studied earlier.

In the canonical model, the risk and its approximation are driven by independent coordinates and by a coordinatewise Stein identity. Under correlation, the coordinates of $\pmb{Z}$ interact through $\pmb{\Omega}$, so both the shrinkage-thresholding design and the risk approximation become coupled across coordinates. A direct reuse of the canonical GEST rule on $\pmb{Z}$ ignores this dependence and can distort the effective shrinkage and thresholding levels. This motivates 
extensions that retain the GEST form while incorporating the geometry induced by $\pmb{\Omega}$. 

\paragraph*{GEST under correlation}
We consider the correlated GEST form
\[
\hat{\pmb{\theta}}_{\mathrm{GEST}}^{+}
=
\left(1-\frac{t(|\pmb{Z}|)}{|\pmb{Z}|^{c(|\pmb{Z}|)}}\right)_+ \circ \hat{\pmb{\theta}},
\]
where $t(|\pmb{Z}|)\ge 0$ and $c(|\pmb{Z}|)\ge 0$, and where $\hat{\pmb{\theta}}$ is a baseline estimator that already reflects the precision structure. We consider two choices for $\hat{\pmb{\theta}}$, namely the MLE and a conditional MLE. 

\subsection{MLE-based shrinkage-thresholding construction}\label{sec:mle_based}

\paragraph*{MLE under invariant quadratic loss}
The maximum likelihood estimator under the model above is
$
\hat{\pmb{\theta}}_{\mathrm{MLE}}
=
\arg\min_{\hat{\pmb{\theta}}}(\pmb{Z}-\hat{\pmb{\theta}})^\top \pmb{\Omega}(\pmb{Z}-\hat{\pmb{\theta}})
=
\pmb{Z}.
$
Applying GEST to the MLE yields the MLE-based shrinkage-thresholding estimator
\[
\hat{\pmb{\theta}}_{\mathrm{GEST,MLE}}
=
\left(1-\frac{t(|\pmb{Z}|)}{|\pmb{Z}|^{c(|\pmb{Z}|)}}\right)\circ \pmb{Z},
\]
and its truncated version is obtained by replacing the multiplier by its positive part.

\paragraph*{Approximate risk in the MLE-based construction}
Given empirical observations $\check{\pmb{z}}$, the invariant risk of $\hat{\pmb{\theta}}_{\mathrm{GEST,MLE}}$ admits a nonparametric approximation analogous to the canonical case. Let
$\tilde{\pmb{z}}=\pmb{\Omega}^{1/2}\check{\pmb{z}}$,
and let $f(\tilde z_i)$ denote the marginal density of $\tilde z_i$ with score $l'(\tilde z_i)$.
Then the approximate risk can be written in the form
\[
\hat{R}(\hat{\pmb{\theta}}_{\mathrm{GEST,MLE}},\pmb{\theta})
=
d+\int_{\mathbb{R}^d}
\hat{H}\!\left(\pmb{z},t,c,\frac{\partial t}{\partial |\pmb{z}|},\frac{\partial c}{\partial |\pmb{z}|}\right)\,d\pmb{z},
\]
where the integrand $\hat{H}$ depends on the correlated geometry through $\pmb{\Omega}$ and on the nonparametric score evaluations through $l'(\tilde{\pmb{z}})$.

\paragraph*{Optimality conditions}
The approximate risk is minimized if and only if the following conditions hold
\begin{equation}
t(|\check{\pmb{z}}|)
=
-\frac{(\sgn(\pmb{z})\circ |\pmb{z}|^{1-c(|\pmb{z}|)})^\top \pmb{\Omega}(\pmb{z}-\tilde{\pmb{z}}-l'(\tilde{\pmb{z}}))}
{(\sgn(\pmb{z})\circ |\pmb{z}|^{1-c(|\pmb{z}|)})^\top \pmb{\Omega}(\sgn(\pmb{z})\circ |\pmb{z}|^{1-c(|\pmb{z}|)})}
\label{eq:approxRiskCond2_1}
\end{equation}
and
\begin{eqnarray}
&&\frac{(\sgn(\pmb{z})\circ |\pmb{z}|^{1-c(|\pmb{z}|)})^\top \pmb{\Omega}(\pmb{z}-\tilde{\pmb{z}}-l'(\tilde{\pmb{z}}))}
{(\sgn(\pmb{z})\circ |\pmb{z}|^{1-c(|\pmb{z}|)})^\top \pmb{\Omega}(\sgn(\pmb{z})\circ |\pmb{z}|^{1-c(|\pmb{z}|)})}  \label{eq:approxRiskCond2_2}\\
&=&
\frac{(\sgn(\pmb{z})\circ |\pmb{z}|^{1-c(|\pmb{z}|)}\circ \ln|\pmb{z}|)^\top \pmb{\Omega}(\pmb{z}-\tilde{\pmb{z}}-l'(\tilde{\pmb{z}}))}
{(\sgn(\pmb{z})\circ |\pmb{z}|^{1-c(|\pmb{z}|)}\circ \ln|\pmb{z}|)^\top \pmb{\Omega}(\sgn(\pmb{z})\circ |\pmb{z}|^{1-c(|\pmb{z}|)})}.\nonumber
\end{eqnarray}
These equations are the correlated analogues of the canonical optimality conditions. They show that the precision matrix weights both the signal-alignment term and the effective scale term, so correlation directly affects the selected shrinkage-thresholding profile.

\paragraph*{Limitations}
The MLE-based construction is 
simple and keeps the estimator in a direct GEST form applied to $\pmb{Z}$. Its limitation is that correlation enters globally through $\pmb{\Omega}$, so strong dependence or ill-conditioning can make the resulting shrinkage-thresholding profile sensitive to covariance misspecification and to marginal score estimation on the transformed scale.

\subsection{Conditional MLE-based construction}\label{sec:cmle_based}

\paragraph*{Conditional MLE}
An alternative to direct shrinkage-thresholding on the MLE is to define the baseline estimator through coordinatewise conditional likelihoods. For each coordinate,
\[
\theta_{\mathrm{CMLE},i}\mid \hat{\pmb{\theta}}_{-i}
=
\arg\min_{\hat{\theta}_i}
(\pmb{Z}-\hat{\pmb{\theta}})^\top \pmb{\Omega}(\pmb{Z}-\hat{\pmb{\theta}})
=
Z_i+\Omega_{ii}^{-1}\pmb{\Omega}_{i,-i}(\pmb{Z}_{-i}-\hat{\pmb{\theta}}_{-i}).
\]
In matrix form,
\[
\hat{\pmb{\theta}}_{\mathrm{CMLE}}
=
\pmb{Z}+\diag(\pmb{\Omega})^{-1}(\pmb{\Omega}-\diag(\Omega))(\pmb{Z}-\hat{\pmb{\theta}}),
\]
or equivalently
\[
\hat{\pmb{\theta}}_{\mathrm{CMLE}}
=
\diag(\pmb{\Omega})^{-1}\pmb{\Omega}\pmb{Z}
-
\diag(\pmb{\Omega})^{-1}(\pmb{\Omega}-\diag(\Omega))\hat{\pmb{\theta}}.
\]
Thus the conditional MLE already incorporates coordinate-specific adjustment through the precision matrix, which makes it a natural baseline for correlated shrinkage-thresholding \cite{besag1974spatial}.

\paragraph*{Restricted CMLE subclass}
In the CMLE-based construction, we restrict attention to the subclass with constant $t$ and $c$. This restriction is not only computationally convenient. It also preserves a tractable fixed-point structure and yields an explicit Stein unbiased risk estimate under the invariant quadratic loss \cite{eldar2009generalized}. If $t(\cdot)$ and $c(\cdot)$ were allowed to vary as general functionals of $|\pmb{Z}|$, then the coupled CMLE system would acquire additional derivative terms through the dependence of $t$ and $c$ on the data, which would substantially complicate both the divergence calculation and the SURE formula. Moreover, the CMLE baseline already induces coordinate-specific adaptation through $\pmb{\Omega}$, so constant $(t,c)$ still defines a flexible correlated shrinkage-thresholding rule rather than a purely global one.

\paragraph*{CMLE-based GEST rule}
Given constants $t\ge 0$ and $c\ge 0$, define $\hat{\pmb{\theta}}_{\mathrm{GEST,CMLE}}^{+}$ and $\hat{\pmb{\theta}}_{\mathrm{CMLE}}$ jointly through
\begin{eqnarray}
&&\hat{\pmb{\theta}}_{\mathrm{GEST,CMLE}}^{+}
=
\left(1-\frac{t}{|\hat{\pmb{\theta}}_{\mathrm{CMLE}}|^c}\right)_+ \circ \hat{\pmb{\theta}}_{\mathrm{CMLE}}\label{eqn:fixed_system}\\
&&\hat{\pmb{\theta}}_{\mathrm{CMLE}}
=
\pmb{Z}
+
\diag(\pmb{\Omega})^{-1}(\pmb{\Omega}-\diag(\Omega))
(\pmb{Z}-\hat{\pmb{\theta}}_{\mathrm{GEST,CMLE}}^{+}).\nonumber
\end{eqnarray}
This defines a fixed-point system in which the conditional MLE supplies the correlation-adjusted baseline and the power-thresholding operator supplies shrinkage and thresholding. In practice, the system is solved iteratively and converges rapidly under the regularity conditions used in our numerical studies.

\paragraph*{Why conditionalization improves stability}
The conditional construction localizes the effect of correlation. Instead of thresholding raw coordinates that may be inflated by dependence, it thresholds conditionally adjusted coordinates whose magnitudes are interpreted relative to the local regression structure encoded by $\pmb{\Omega}$. This reduces instability caused by strong off-diagonal dependence and yields a more interpretable thresholding mechanism when the precision matrix is sparse or approximately sparse \cite{besag1974spatial,meinshausen2006high}.

\paragraph*{SURE for the CMLE-based construction}
A key advantage of the restricted CMLE subclass is that it admits an explicit SURE under the invariant quadratic risk \cite{eldar2009generalized}. For $\hat{\pmb{\theta}}_{\mathrm{GEST,CMLE}}^{+}$,
\[
\nabla\cdot \hat{\pmb{\theta}}_{\mathrm{GEST,CMLE}}^{+}
=
\tr\!\left[
\left(\pmb{I}_d+\pmb{\Psi}\diag(\pmb{\Omega})^{-1}-\pmb{\Psi}\right)^{-1}
\pmb{\Psi}\diag(\pmb{\Omega})^{-1}\pmb{\Omega}
\right],
\]
where
\[
\pmb{\Psi}
=
\diag\!\left[
\left(
\frac{\hat{\pmb{\theta}}_{\mathrm{GEST,CMLE}}^{+}}{\hat{\pmb{\theta}}_{\mathrm{CMLE}}}
+
\frac{tc|\hat{\pmb{\theta}}_{\mathrm{CMLE}}|^{-c}}{\diag(\pmb{\Omega})^{c/2}}
\right)
1_{\hat{\pmb{\theta}}_{\mathrm{GEST,CMLE}}^{+}\neq 0}
\right].
\]
The corresponding Stein unbiased risk estimate is
\[
\mathrm{SURE}(\hat{\pmb{\theta}}_{\mathrm{GEST,CMLE}}^{+})
=
-d
+
2\nabla\cdot \hat{\pmb{\theta}}_{\mathrm{GEST,CMLE}}^{+}
+
(\hat{\pmb{\theta}}_{\mathrm{GEST,CMLE}}^{+}-\pmb{Z})^\top
\pmb{\Omega}
(\hat{\pmb{\theta}}_{\mathrm{GEST,CMLE}}^{+}-\pmb{Z}),
\]
and it satisfies
\[
R(\hat{\pmb{\theta}},\pmb{\theta})
=
E_{\pmb{\theta}}\!\left[\mathrm{SURE}(\hat{\pmb{\theta}})\right].
\]
This provides a direct approximate-risk criterion for selecting $(t,c)$ in the CMLE-based construction. It also explains why the constant-parameter subclass is especially attractive in the correlated setting, since it preserves an explicit and computationally usable risk estimate.

\subsection{Theoretical results in the correlated setting}\label{sec:correlated_theory}

We now state the main theoretical properties for the correlated constructions. The proof are deferred to the Supplement \cite{jiang2026nomadS}. 

\paragraph*{Admissible classes}
For the MLE-based construction, $\mathcal{F}_{\mathrm{MLE}}$ is a compact class of measurable bounded shrinkage-thresholding functionals $(t,c)$ of $|\pmb Z|$ such that $c\in[0,1]$, the partial derivatives of $t$ and $c$ exist almost everywhere and have at most polynomial growth, the correlated approximate-risk integrand is integrable, and the transformed marginal score is well defined with at most polynomial growth.
Thus, in the MLE-based construction, both $t$ and $c$ are allowed to depend on the observed magnitude profile $|\pmb{Z}|$, subject to boundedness and the smoothness conditions needed for the approximate-risk criterion in \eqref{eq:approxRiskCond2_1}--\eqref{eq:approxRiskCond2_2}. For the CMLE-based construction, we restrict attention to the constant-parameter subclass
\[
\mathcal{F}_{\mathrm{CMLE}}=[0,t_{\max}]\times[0,1].
\]

\paragraph*{Standing assumptions}
Assume that $\pmb{\Omega}$ is symmetric positive definite and its eigenvalues are uniformly bounded away from zero and infinity, namely
\[
0<\lambda_{\min}\le \lambda_{\min}(\pmb{\Omega})\le \lambda_{\max}(\pmb{\Omega})\le \lambda_{\max}<\infty
\]
for all $d$. Assume also that
\[
\frac{1}{d}\|\pmb{\Omega}^{1/2}\pmb{\theta}\|_2^2\le B_2
\]
for a constant $B_2<\infty$.

For the MLE-based construction, assume the score estimator used in \eqref{eq:approxRiskCond2_1}--\eqref{eq:approxRiskCond2_2} is uniformly consistent on compact sets for the transformed coordinates $\tilde{\pmb{z}}=\pmb{\Omega}^{1/2}\check{\pmb{z}}$, and assume the sieve approximation and uniform approximation conditions from Section~\ref{sec:canonical_theory} hold after replacing the canonical quadratic loss by the invariant quadratic loss.

For the CMLE-based construction, assume
\[
0<\omega_{\min}\le \Omega_{ii}\le \omega_{\max}<\infty
\qquad\text{for all } i,
\]
and
\[
\max_{1\le i\le d}\sum_{j\neq i}\left|\frac{\Omega_{ij}}{\Omega_{ii}}\right|\le \rho<1.
\]
Then the CMLE fixed-point map is a contraction uniformly over $(t,c)\in\mathcal{F}_{\mathrm{CMLE}}$. Write
\[
\pmb{A}=\diag(\pmb\Omega)^{-1}\bigl(\pmb{\Omega}-\diag(\pmb\Omega)\bigr).
\]
and assume \(\|\pmb A\|_{\mathrm{op}}\le \rho_B<1\) for all sufficiently large \(d\). In addition, for every \((t,c)\in\mathcal F_{\mathrm{CMLE}}\), the matrix
\[
\pmb M(t,c):=\pmb I_d+\pmb\Psi(t,c)\diag(\pmb\Omega)^{-1}-\pmb\Psi(t,c)
\]
is uniformly invertible in the sense that there exists \(\delta_0>0\) such that
\[
\inf_{(t,c)\in\mathcal F_{\mathrm{CMLE}}}
\lambda_{\min}\!\bigl(\pmb M(t,c)\bigr)\ge \delta_0
\qquad\text{a.s.}
\]

\paragraph*{Consistency of Data-driven MLE-based NOMAD}
For $(t,c)\in\mathcal{F}_{\mathrm{MLE}}$, let $\hat{\pmb{\theta}}_{t,c}^{\mathrm{MLE}}$ denote the corresponding MLE-based correlated GEST estimator and define
\[
R_{\mathrm{MLE}}(t,c;\pmb{\theta},\pmb{\Omega})
=
E_{\pmb{\theta}}
\left[
(\hat{\pmb{\theta}}_{t,c}^{\mathrm{MLE}}-\pmb{\theta})^\top
\pmb{\Omega}
(\hat{\pmb{\theta}}_{t,c}^{\mathrm{MLE}}-\pmb{\theta})
\right].
\]
Let
$
(\hat t_{\mathrm{MLE}},\hat c_{\mathrm{MLE}})
\in
\arg\min_{(t,c)\in\mathcal{F}_{\mathrm{MLE}}}
\hat R_{\mathrm{MLE}}(t,c;\check{\pmb{z}},\pmb{\Omega}),
$
with a measurable tie-breaking rule if needed, and define
$
\hat{\pmb{\theta}}_{\mathrm{NOMAD,MLE}}
=
\hat{\pmb{\theta}}_{\hat t_{\mathrm{MLE}},\hat c_{\mathrm{MLE}}}^{\mathrm{MLE}}.
$
Also define the oracle risk over the MLE-based class by
\[
R_{\mathrm{MLE}}^\star(\pmb{\theta},\pmb{\Omega})
=
\inf_{(t,c)\in\mathcal{F}_{\mathrm{MLE}}}
R_{\mathrm{MLE}}(t,c;\pmb{\theta},\pmb{\Omega}).
\]

\begin{theorem}[Oracle-risk consistency of the data-driven MLE-based NOMAD estimator]\label{thm:consistency_mle_corr}
Assume the standing assumptions above. Then
\[
R_{\mathrm{MLE}}(\hat t_{\mathrm{MLE}},\hat c_{\mathrm{MLE}};\pmb{\theta},\pmb{\Omega})
-
R_{\mathrm{MLE}}^\star(\pmb{\theta},\pmb{\Omega})
=
o_{\mathbb{P}}(d).
\]
\end{theorem}

\paragraph*{Consistency of Data-driven CMLE-based NOMAD}
For each $(t,c)\in\mathcal{F}_{\mathrm{CMLE}}$, 
define the invariant quadratic risk
\[
R_{\mathrm{CMLE}}(t,c;\pmb{\theta},\pmb{\Omega})
=
E_{\pmb{\theta}}
\left[
\bigl(\hat{\pmb{\theta}}_{\mathrm{GEST,CMLE}}^{+}(t,c)-\pmb{\theta}\bigr)^\top
\pmb{\Omega}
\bigl(\hat{\pmb{\theta}}_{\mathrm{GEST,CMLE}}^{+}(t,c)-\pmb{\theta}\bigr)
\right],
\]
and the oracle risk over the CMLE-based class by
\[
R_{\mathrm{CMLE}}^\star(\pmb{\theta},\pmb{\Omega})
=
\inf_{(t,c)\in\mathcal{F}_{\mathrm{CMLE}}}
R_{\mathrm{CMLE}}(t,c;\pmb{\theta},\pmb{\Omega}).
\]

Define the SURE criterion
\[
\mathrm{SURE}(t,c):=
\mathrm{SURE}\bigl(\hat{\pmb{\theta}}_{\mathrm{GEST,CMLE}}^{+}(t,c)\bigr),
\]
and let
$(\hat t_{\mathrm{CMLE}},\hat c_{\mathrm{CMLE}})
\in
\arg\min_{(t,c)\in\mathcal{F}_{\mathrm{CMLE}}}
\mathrm{SURE}(t,c)$, 
with a measurable tie-breaking rule if needed. The data-driven CMLE-based NOMAD estimator is
$
\hat{\pmb{\theta}}_{\mathrm{NOMAD,CMLE}}
=
\hat{\pmb{\theta}}_{\mathrm{GEST,CMLE}}^{+}(\hat t_{\mathrm{CMLE}},\hat c_{\mathrm{CMLE}}).
$

\begin{theorem}[Oracle-risk consistency of the data-driven CMLE-based NOMAD estimator]\label{thm:consistency_cmle_corr}
Assume the standing assumptions above. Then
\[
R_{\mathrm{CMLE}}(\hat t_{\mathrm{CMLE}},\hat c_{\mathrm{CMLE}};\pmb{\theta},\pmb{\Omega})
-
R_{\mathrm{CMLE}}^\star(\pmb{\theta},\pmb{\Omega})
=
o_{\mathbb{P}}(d).
\]
\end{theorem}

\paragraph*{Approximate-risk inequalities}
Each construction satisfies a deterministic ordering property on its own criterion scale, which follows directly from the definition of the corresponding estimator as a minimizer of its empirical criterion over the stated admissible class. 

\begin{theorem}[Approximate-risk inequalities]\label{thm:risk_ineq_corr}
Let $(\hat t_{\mathrm{MLE}},\hat c_{\mathrm{MLE}})$ and $(\hat t_{\mathrm{CMLE}},\hat c_{\mathrm{CMLE}})$ be defined as above. Then
\[
\hat{R}_{\mathrm{MLE}}(\hat t_{\mathrm{MLE}},\hat c_{\mathrm{MLE}};\check{\pmb{z}},\pmb{\Omega})
\le
\hat{R}_{\mathrm{MLE}}(t,c;\check{\pmb{z}},\pmb{\Omega})
\qquad
\text{for all } (t,c)\in\mathcal{F}_{\mathrm{MLE}},
\]
and
\[
\mathrm{SURE}(\hat t_{\mathrm{CMLE}},\hat c_{\mathrm{CMLE}})
\le
\mathrm{SURE}(t,c)
\qquad
\text{for all } (t,c)\in\mathcal{F}_{\mathrm{CMLE}}.
\]
In particular,
\[
\hat{R}_{\mathrm{MLE}}(\hat{\pmb{\theta}}_{\mathrm{NOMAD,MLE}},\pmb{\theta})
\le
\hat{R}_{\mathrm{MLE}}(\hat{\pmb{\theta}}_{\mathrm{JS}},\pmb{\theta})
\le
\hat{R}_{\mathrm{MLE}}(\hat{\pmb{\theta}}_{\mathrm{MLE}},\pmb{\theta})
=
d,
\]
and $\mathrm{SURE}(\hat t_{\mathrm{CMLE}},\hat c_{\mathrm{CMLE}})
\le
\mathrm{SURE}(t, 1)$ which corresponds to lasso.
\end{theorem}

\paragraph*{Broad structural assumptions imply asymptotic equivalence on a common class}
To compare the MLE-based and CMLE-based constructions on a common scale, restrict both procedures to the same constant-parameter class
\[
\mathcal F_{\mathrm{CMLE}}=[0,t_{\max}]\times[0,1].
\]
Impose the following additional structural assumptions:
\begin{enumerate}
\item the conditional coupling matrix satisfies
$
\|\pmb{A}\|_{\mathrm{op}}=o(1)
$
\item the weak-correction condition holds:
\begin{equation}
E_{\pmb{\theta}}\!\left[\|\pmb{A}(\pmb{Z}-\pmb{\theta})\|_2^2\right]=o(d)
\label{eq:weak_correction_condition_section8}
\end{equation}
\end{enumerate}
These conditions describe a weak-dependence regime in which the CMLE correction is stable and its stochastic contribution is lower order on the scale of \(d\).

For each \((t,c)\in\mathcal{F}_{\mathrm{CMLE}}\), let
$
\hat{\pmb{\theta}}_{t,c}^{\mathrm{MLE}}
=
\mathcal{P}_{t,c}(\pmb{Z})$,
$\hat{\pmb{\theta}}_{t,c}^{\mathrm{CMLE}}
=
\hat{\pmb{\theta}}_{\mathrm{GEST,CMLE}}^{+}(t,c).
$
Let
\[
(\hat t_{\mathrm{MLE},0},\hat c_{\mathrm{MLE},0})
\in
\arg\min_{(t,c)\in\mathcal F_{\mathrm{CMLE}}}
\hat R_{\mathrm{MLE}}(t,c;\check{\pmb z},\pmb\Omega),
\]
and define
$
\hat{\pmb{\theta}}_{\mathrm{NOMAD,MLE},0}
=
\hat{\pmb{\theta}}_{\hat t_{\mathrm{MLE},0},\hat c_{\mathrm{MLE},0}}^{\mathrm{MLE}}.
$

\begin{theorem}[Asymptotic equivalence of the CMLE- and MLE-based constructions on the common class]
\label{thm:cmle_mle_equiv}
Assume the standing assumptions and broad structural assumptions listed above. Then 
\[
\sup_{(t,c)\in\mathcal{F}_{\mathrm{CMLE}}}
\left|
R_{\mathrm{CMLE}}(t,c;\pmb{\theta},\pmb{\Omega})
-
R_{\mathrm{MLE}}(t,c;\pmb{\theta},\pmb{\Omega})
\right|
=
o(d).
\]
Consequently,
\[
\left|
R_{\mathrm{CMLE}}(\hat{\pmb{\theta}}_{\mathrm{NOMAD,CMLE}},\pmb{\theta})
-
R_{\mathrm{MLE}}(\hat{\pmb{\theta}}_{\mathrm{NOMAD,MLE},0},\pmb{\theta})
\right|
=
o_{\mathbb{P}}(d).
\]
\end{theorem}

\begin{remark}
Compared with the canonical setting, the correlated theory is naturally stated in terms of $\pmb{\Omega}$-weighted oracle-risk consistency rather than unweighted Euclidean consistency. Theorems~\ref{thm:consistency_mle_corr} and \ref{thm:consistency_cmle_corr} show that both correlated NOMAD constructions are oracle-risk consistent relative to their own admissible shrinkage-thresholding classes. Theorem~\ref{thm:cmle_mle_equiv} shows that under broad structural assumptions, namely weak dependence, stable conditional correction, and sparse signal, the CMLE-based and MLE-based NOMAD estimators are asymptotically equivalent in invariant quadratic risk. 
The MLE-based construction parallels the canonical framework more closely, which is simpler and is convenient when $\pmb{\Omega}$ is well-conditioned and accurately known or estimated. In contrast, the CMLE-based construction restricts to a constant-parameter subclass in order to preserve an explicit SURE formula and a tractable fixed-point representation. This tradeoff yields a stable and practically useful correlated extension and prepares the regression formulation in the next section. In regimes with structured dependence, the conditional adjustment may reduce spurious cross-talk into null coordinates before thresholding. Establishing general strict dominance requires additional directional assumptions and is beyond the scope of the paper. 
\end{remark}

\section{Linear Regression and Penalized Regression Representation}\label{sec:regression}

This section connects the correlated normal mean extension to linear regression. The regression model induces a multivariate normal mean problem for suitable sufficient statistics, and the shrinkage-thresholding map selected by approximate risk minimization yields a regression estimator. We then show that the resulting estimator admits an equivalent penalized regression formulation with an explicit data-adaptive penalty induced by the inverse of the shrinkage-thresholding operator \cite{antoniadis2001regularization}.

\subsection{Linear regression as an induced normal mean problem}\label{sec:regression_model}

\paragraph*{Model setup}
Consider the Gaussian linear model
\begin{equation}
\pmb{y}=\pmb{X}\pmb{\beta}+\pmb{\varepsilon},
\qquad
\pmb{\varepsilon}\sim N(\pmb{0},\sigma^2\pmb{I}_n),
\label{eq:linreg_model}
\end{equation}
where $\pmb{y}\in\mathbb{R}^n$ is the response, $\pmb{X}\in\mathbb{R}^{n\times p}$ is the design matrix, $\pmb{\beta}\in\mathbb{R}^p$ is the unknown coefficient vector, and $\sigma^2>0$ is the noise variance.

\paragraph*{General design reduction to a correlated normal mean model}
For a general full column rank design, the least squares estimator satisfies
\[
\hat{\pmb{\beta}}_{\mathrm{LS}}=(\pmb{X}^\top\pmb{X})^{-1}\pmb{X}^\top\pmb{y},
\qquad
\hat{\pmb{\beta}}_{\mathrm{LS}}\sim N\!\left(\pmb{\beta},\sigma^2(\pmb{X}^\top\pmb{X})^{-1}\right).
\]
Define the standardized statistic
$\pmb{Z}=\frac{1}{\sigma}\hat{\pmb{\beta}}_{\mathrm{LS}}$.
Then
\begin{equation}
\pmb{Z}\sim N(\pmb{\theta},\pmb{\Sigma}),
\qquad
\pmb{\theta}=\frac{1}{\sigma}\pmb{\beta},
\qquad
\pmb{\Sigma}=(\pmb{X}^\top\pmb{X})^{-1}
\label{eq:reg_as_correlated_means}
\end{equation}
This representation links regression to the correlated multivariate normal mean model studied in Section~\ref{sec:correlated_extension}.

\subsection{NOMAD estimator for linear regression}\label{sec:regression_nomad}

\begin{figure}[h!]
    \centering
    \begin{subfigure}[b]{0.35\textwidth}
        \centering
        \includegraphics[width=\textwidth]{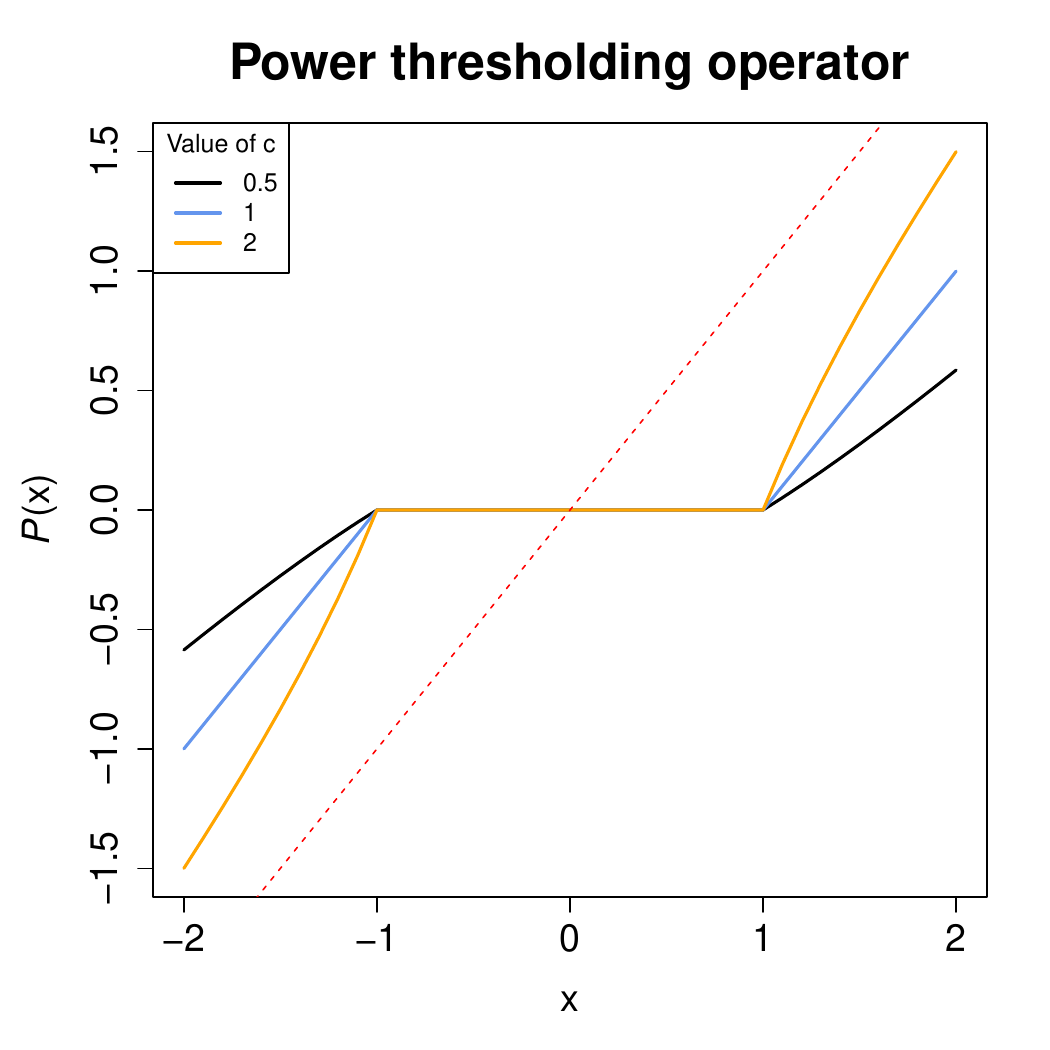}
    \end{subfigure}
    \begin{subfigure}[b]{0.35\textwidth}
        \centering
        \includegraphics[width=\textwidth]{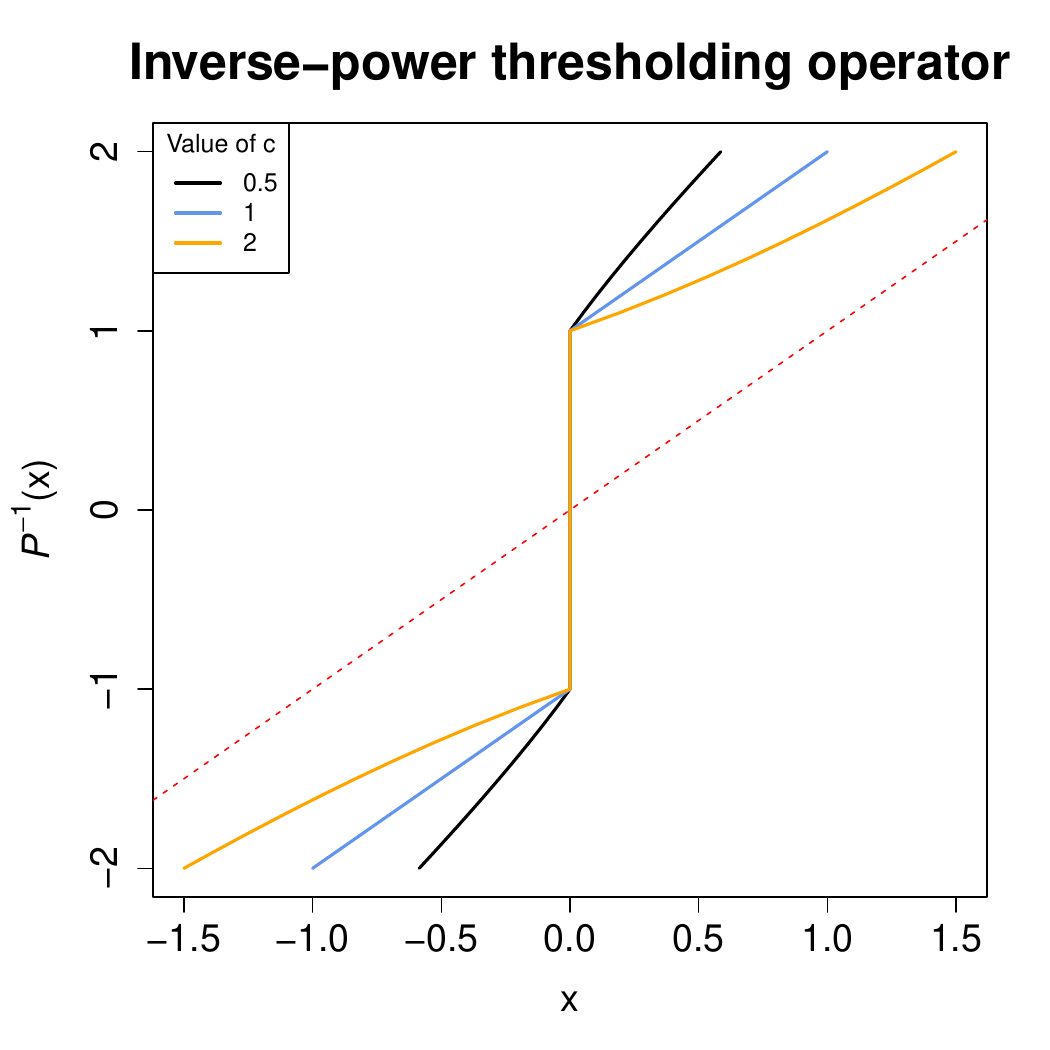}
    \end{subfigure}
    
    \caption{Power Thresholding operator $\mathcal{P}_{t,c}(x)$ (left) and Inverse Power Thresholding Operator $\mathcal{P}^{-1}_{t,c}$ (right) when $t=1$.}
    \label{fig:power_thresholding_operator}
\end{figure}

\paragraph*{Power thresholding operator}
For a scalar variable $x$, define the power thresholding operator
\begin{equation}
\mathcal{P}_{t,c}(x)
=
\left(1-\frac{t}{|x|^c}\right)_+x,
\qquad t\ge 0,\ c\ge 0
\label{eq:power_threshold_op}
\end{equation}
This operator matches the coordinatewise form of the truncated GEST rule and it is monotone nondecreasing in $x$ for fixed $(t,c)$, so it admits an inverse map on its range, which we denote by $\mathcal{P}^{-1}_{t,c}$. Both operators with $t=1$ and different values of $c$ can be visualized in Figure \ref{fig:power_thresholding_operator}.

\paragraph*{Regression NOMAD}
Let $\check{\pmb{z}}=\hat{\pmb{\beta}}_{\mathrm{LS}}/\hat{\sigma}$, where $\hat{\sigma}$ is a noise estimator, and let $\pmb{\Sigma}=(\pmb{X}^\top\pmb{X})^{-1}$ and $\pmb{\Omega}=\pmb{\Sigma}^{-1}=\pmb{X}^\top\pmb{X}$. CMLE-based procedure is more stable and it typically yields lower risk presented in Section~\ref{sec:sims_correlated}. We apply the CMLE-based correlated normal mean construction of Section~\ref{sec:correlated_extension} to \eqref{eq:reg_as_correlated_means} and obtain an approximate-risk minimizer $(\tilde t,\tilde c)$, which defines a coordinatewise shrinkage-thresholding map through $\mathcal{P}_{\tilde t,\tilde c}$. The resulting regression estimator is
$\hat{\pmb{\beta}}_{\mathrm{NOMAD}}=\hat{\sigma}\,\hat{\pmb{\theta}}_{\mathrm{NOMAD}}$,
where $\hat{\pmb{\theta}}_{\mathrm{NOMAD}}$ is the estimated mean vector in the induced normal mean model computed by the CMLE-based procedure in Section~\ref{sec:correlated_extension}.

\subsection{Equivalent penalized regression formulation}\label{sec:regression_penalized}

The shrinkage-thresholding map $\mathcal{P}_{t,c}$ induces a coordinate-separable penalty through its inverse. Define the scalar penalty function
\begin{equation}
\rho_{t,c}(v)
=
\int_{0}^{v}\bigl(\mathcal{P}^{-1}_{t,c}(u)-u\bigr)\,du
\label{eq:rho_def}
\end{equation}
Whenever $\mathcal{P}_{t,c}$ is monotone, $\mathcal{P}^{-1}_{t,c}$ is well-defined on the range and $\rho_{t,c}$ is nonnegative with $\rho_{t,c}(0)=0$ \cite{antoniadis2001regularization}.

The next result gives an explicit penalized regression objective whose solution coincides with the NOMAD estimator in linear regression. We defer the proof to Supplement \cite{jiang2026nomadS}. 

\begin{theorem}\label{thm:penalized_regression_nomad}
Fix $(t,c)$ and define $\mathcal{P}_{t,c}$ by \eqref{eq:power_threshold_op} and $\rho_{t,c}$ by \eqref{eq:rho_def}. The NOMAD estimator in linear regression is equivalent to solving
\begin{equation}
\min_{\pmb{\beta}\in\mathbb{R}^p}
\ \frac{1}{2\sigma^2}\|\pmb{y}-\pmb{X}\pmb{\beta}\|_2^2
+\sum_{i=1}^p
\int_{0}^{\frac{(\pmb{X}^\top\pmb{X})_{ii}^{1/2}\beta_i}{\sigma}}
\bigl(\mathcal{P}^{-1}_{t,c}(u)-u\bigr)\,du
\label{eq:penalized_nomad}
\end{equation}
If $0\le c\le 1$, the objective in \eqref{eq:penalized_nomad} is convex.
If $c>1$, the induced penalty is concave.
\end{theorem}

\subsection{Connections to ridge and lasso}\label{sec:regression_ridge_lasso}

\paragraph*{Lasso as the case $c=1$}
If $c=1$, then $\mathcal{P}_{t,1}(x)=(1-t/|x|)_+x$ is soft-thresholding and the induced penalty in \eqref{eq:penalized_nomad} reduces to an $\ell_1$ penalty. In particular, \eqref{eq:penalized_nomad} coincides with the lasso objective with tuning parameter $\lambda=t$ \cite{fan2001variable}.

\paragraph*{Ridge-type shrinkage as a shrinkage-only regime}
When $c=0$ and the shrinkage factor does not cross zero, the map becomes linear shrinkage and the induced penalty is quadratic, which yields ridge-type behavior. In this regime NOMAD behaves as a shrinkage-only estimator and the penalized representation reduces to ridge regression after matching the shrinkage level to the quadratic penalty parameter \cite{frank1993statistical,fu1998penalized}.

\paragraph*{Penalty shapes and adaptivity}
For $0<c<1$, the penalty interpolates between quadratic shrinkage and absolute-value thresholding, and it yields shrinkage-thresholding behavior that is weaker than soft-thresholding near the origin and closer to shrinkage on larger magnitudes. For $c>1$, concavity induces nonconvex regularization and it can yield more aggressive sparsity promotion \cite{fu1998penalized,fan2001variable}. The approximate-risk minimization step selects $(t,c)$ from data, so NOMAD adaptively chooses a penalty shape that reflects the empirical distribution of the sufficient statistics.

\section{Simulation Studies}\label{sec:simulations}

This section reports simulation results that assess the finite-sample behavior of NOMAD and examine whether the empirical patterns agree with the theoretical development. We consider three settings. We begin with the canonical normal mean model, where the approximate-risk construction is derived directly. We then study correlated normal mean models, where we compare the MLE-based and CMLE-based extensions. We finally consider linear regression, where NOMAD induces a data-adaptive penalized estimator. All reported values are averages over repeated Monte Carlo replications.

\subsection{Canonical normal mean simulations}\label{sec:sims_canonical}

We generate
$\pmb{Z}\sim N(\pmb{\theta},\pmb{I}_d)$
with dimension $d=500$. The signal vector $\pmb{\theta}$ is generated coordinatewise from the spike-and-slab model
\[
\theta_i \sim (1-\pi)\delta_0+\pi N(0,\sigma^2),
\]
where $\pi$ controls the proportion of nonzero coordinates and $\sigma^2$ controls the signal magnitude \cite{mitchell1988bayesian,george1993variable}. We vary $\pi$ over a grid from $0$ to $1$. We take $\sigma^2\in\{1,2,3\}$, which correspond to weak, moderate, and strong signal regimes.

\begin{figure}[h!]
    \centering
    \begin{subfigure}[b]{0.32\textwidth}
        \centering
        \includegraphics[width=\textwidth]{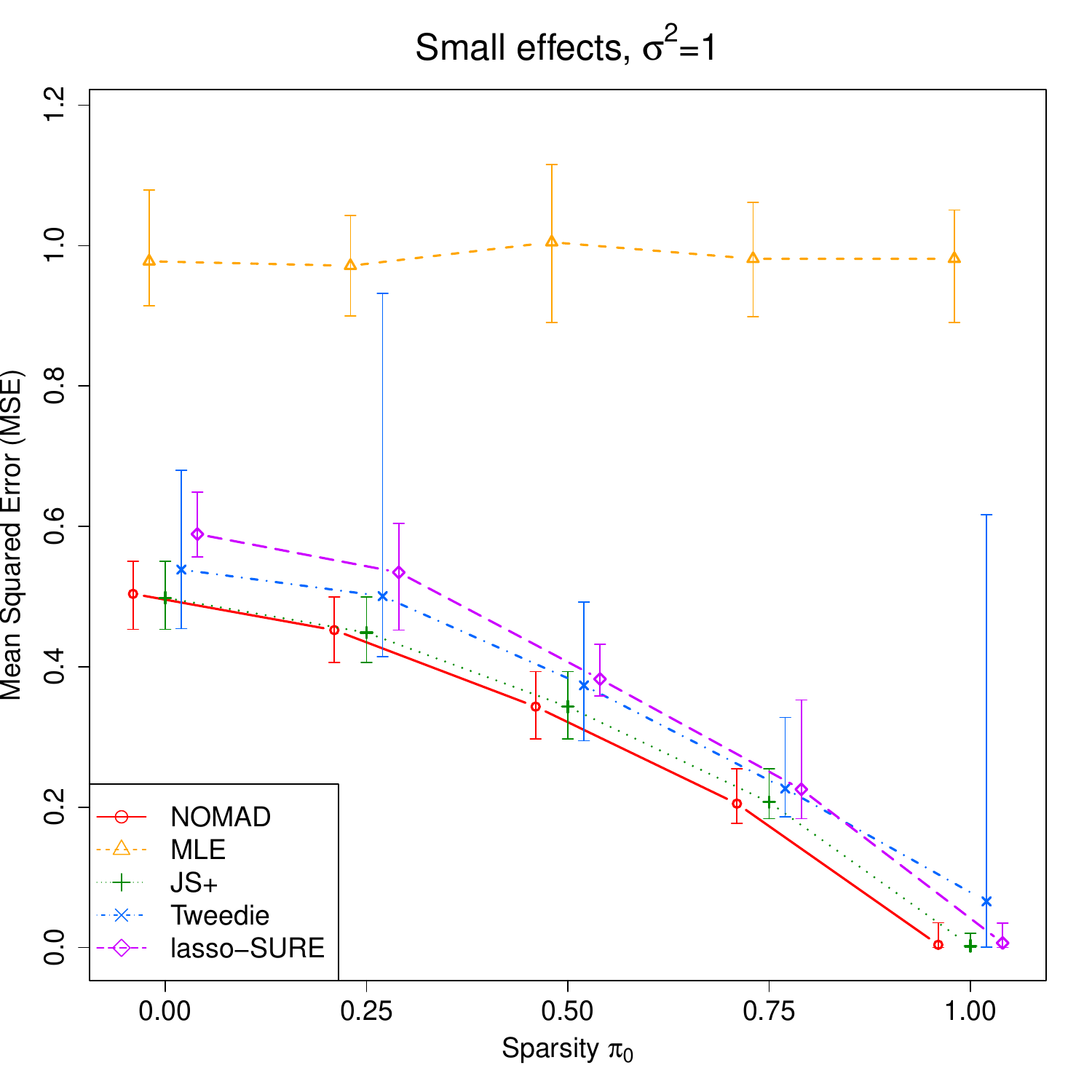}
    \end{subfigure}
    \begin{subfigure}[b]{0.32\textwidth}
        \centering
        \includegraphics[width=\textwidth]{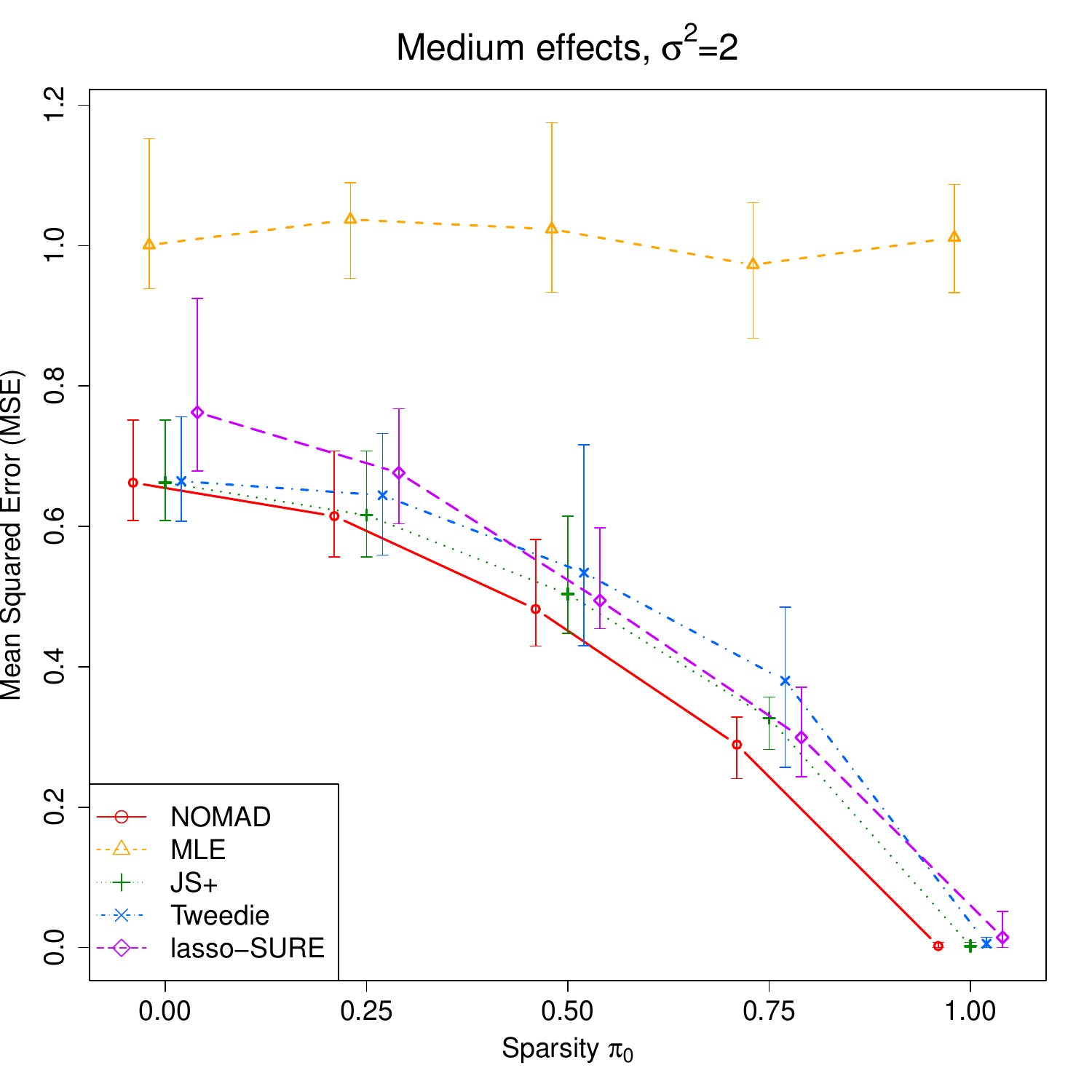}
    \end{subfigure}
    \begin{subfigure}[b]{0.32\textwidth}
        \centering
        \includegraphics[width=\textwidth]{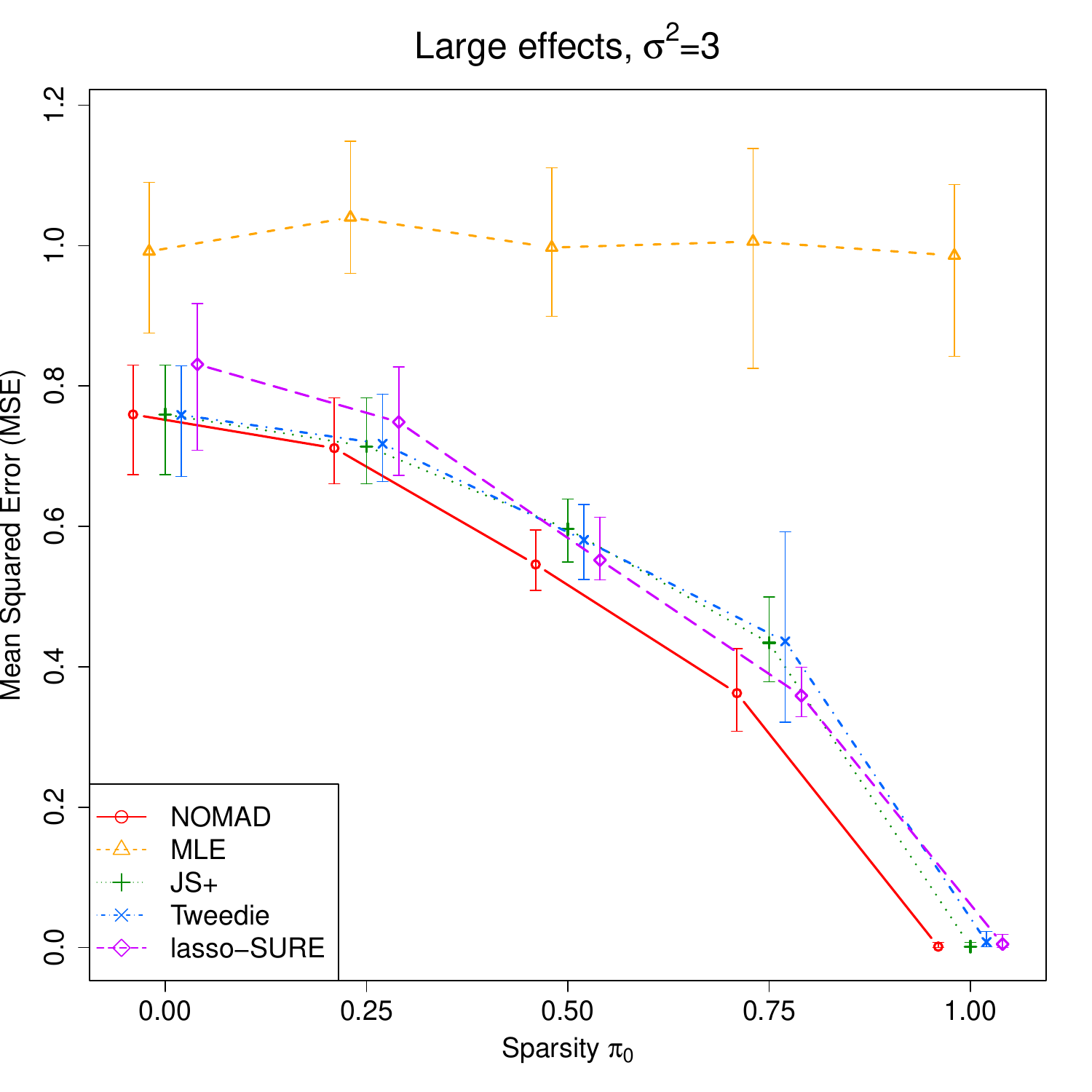}
    \end{subfigure}
    \caption{Per-coordinate mean squared error in the canonical normal mean simulations for weak, moderate, and strong signal regimes.}
    \label{fig:canonical_mse}
\end{figure}

We compare the following procedures:
(a) the maximum likelihood estimator; 
(b) the positive-part James--Stein estimator;
(c) Tweedie's formula based estimator from Lemma~\ref{lem:tweedie}
(d) soft-thresholding with the threshold chosen by SURE;
and (e) NOMAD based on approximate-risk minimization over the GEST class, 
and report the per-coordinate quadratic risk. 
The results are displayed in Figure~\ref{fig:canonical_mse}.

Several patterns are clear. First, the maximum likelihood estimator has the largest risk throughout, with per-coordinate mean squared error close to one, which is consistent with the canonical noise level. Second, soft-thresholding improves substantially in sparse regimes, especially when the fraction of nonzero coordinates is small. Third, the positive-part James--Stein estimator is more competitive in denser regimes, and it often outperforms Tweedie's formula when signals are weak or moderate. In the strong-signal regime, the performances of James--Stein and Tweedie become closer. Across the moderate- and strong-signal regimes, NOMAD yields the smallest risk uniformly over the sparsity range considered. In the weak-signal regime, NOMAD remains competitive and closely tracks the positive-part James--Stein estimator while still improving over the remaining benchmarks. These results agree with the theoretical picture that NOMAD adapts between shrinkage-dominant and thresholding-dominant behavior according to the empirical signal configuration.

\subsection{Correlated multivariate normal mean simulations}\label{sec:sims_correlated}

We generate
$
\pmb{Z}\sim N(\pmb{\theta},\pmb{\Sigma})
$
with dimension $d=500$. We take an AR(1) covariance structure
$
\Sigma_{ij}=\rho^{|i-j|}$,
$\rho=0.5.
$
The signal vector $\pmb{\theta}$ is generated from the same spike-and-slab model as in Section~\ref{sec:sims_canonical}. We again vary the sparsity parameter over a grid from $0$ to $1$ and take $\sigma^2\in\{1,2,3\}$.

\begin{figure}[h!]
    \centering
    \begin{subfigure}[b]{0.32\textwidth}
        \centering
        \includegraphics[width=\textwidth]{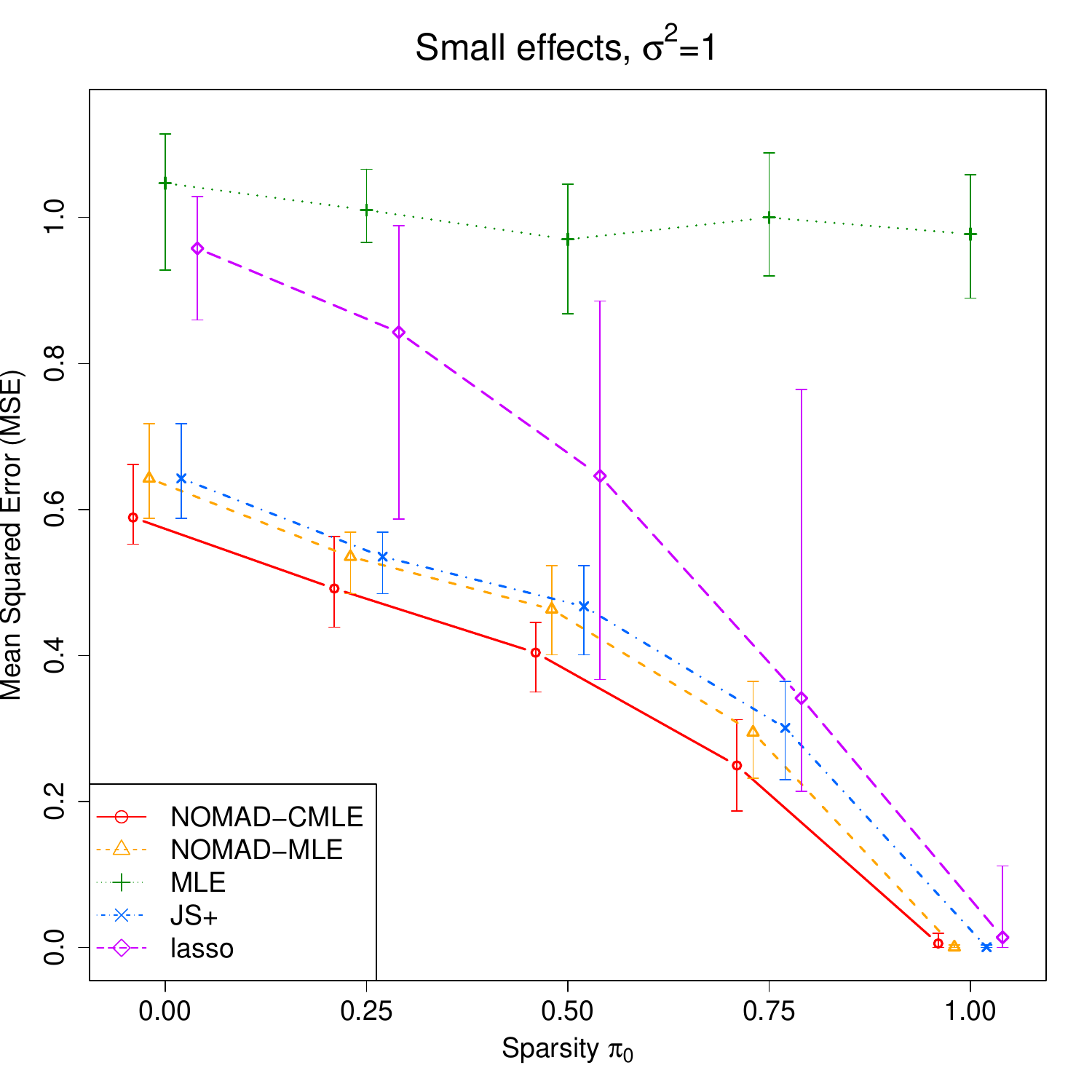}
    \end{subfigure}
    \begin{subfigure}[b]{0.32\textwidth}
        \centering
        \includegraphics[width=\textwidth]{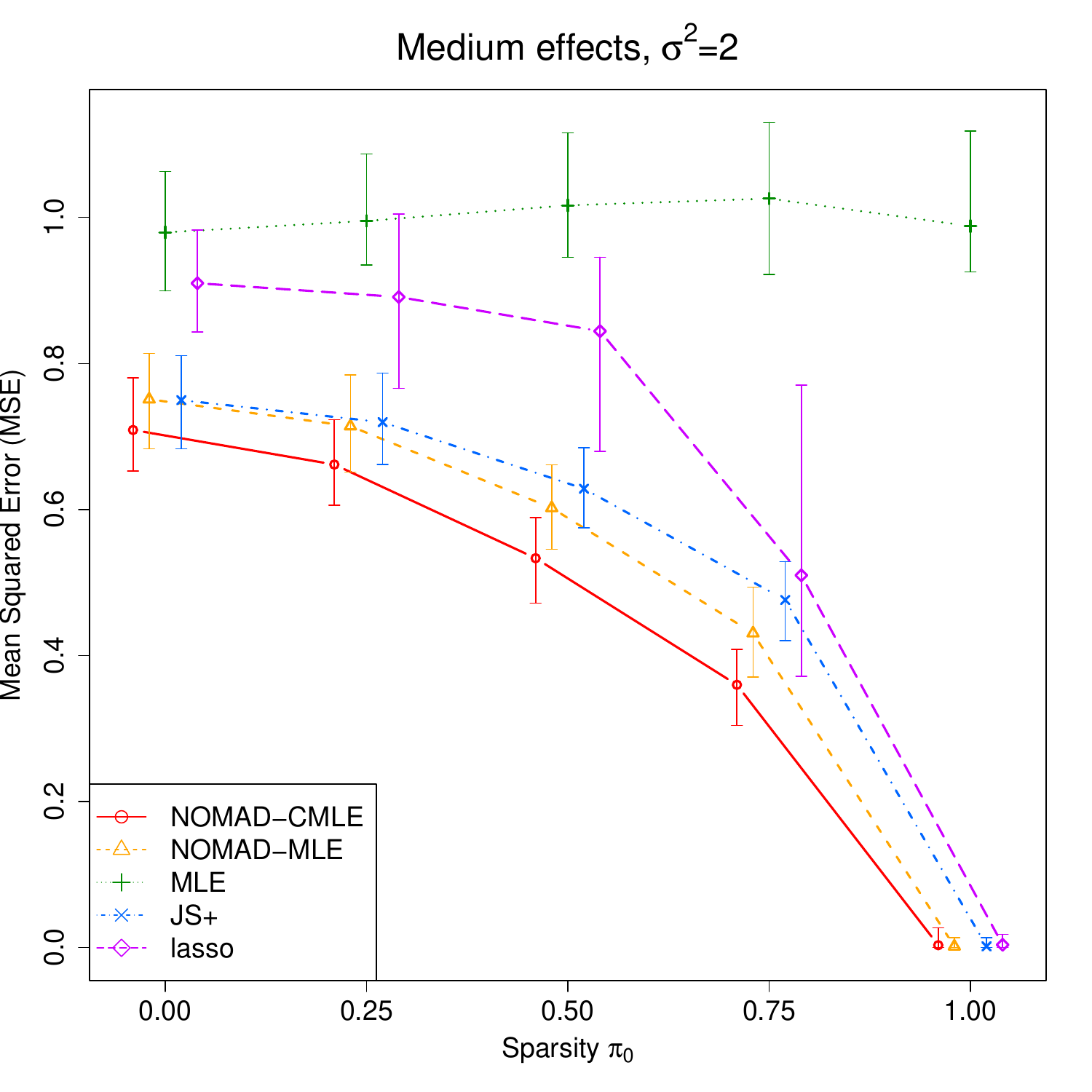}
    \end{subfigure}
    \begin{subfigure}[b]{0.32\textwidth}
        \centering
        \includegraphics[width=\textwidth]{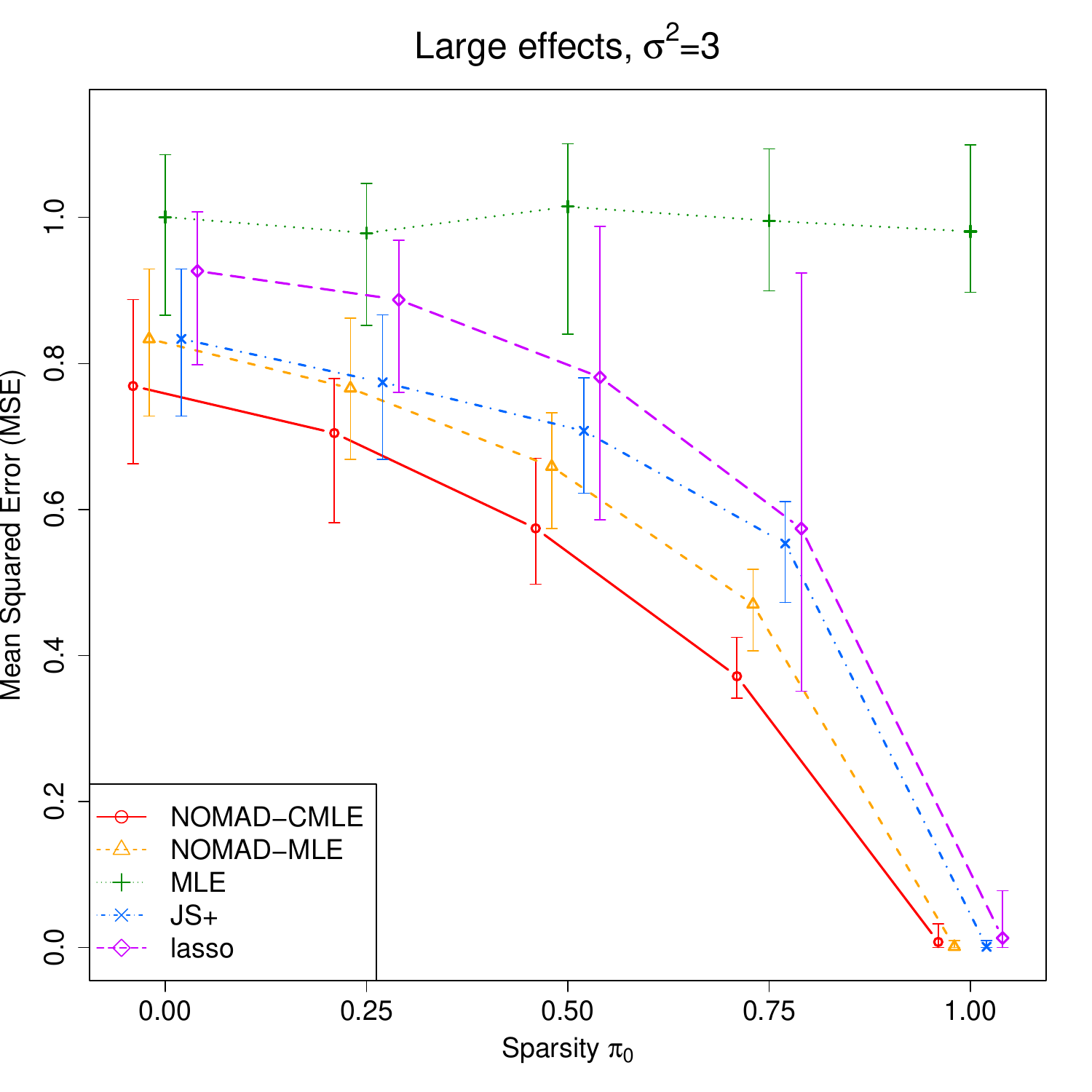}
    \end{subfigure}
    \caption{Per-coordinate mean squared error in the correlated normal mean simulations under AR(1) dependence.}
    \label{fig:ar1_mse}
\end{figure}

We compare the following methods:
(a) maximum likelihood estimator;
(b) whitening-based positive-part James--Stein;
(c) lasso;
(d) NOMAD-MLE, the MLE-based correlated extension;
and (e) NOMAD-CMLE, the CMLE-based correlated extension, 
and report the quadratic risk and its per-coordinate version, estimated by Monte Carlo replication. The results are summarized in Figure~\ref{fig:ar1_mse}.

The maximum likelihood estimator again has the largest risk across all regimes. Under correlation, the whitening-based positive-part James--Stein estimator is generally more competitive than lasso, although the gap narrows as the signal becomes sparser. Both correlated NOMAD procedures improve on the classical benchmarks. Among the two, NOMAD-CMLE is typically more stable and usually attains the smaller risk. This pattern is consistent with the theoretical motivation for conditional adjustment, namely that the CMLE-based construction removes part of the correlation-induced leakage before thresholding. When signals are weak, NOMAD-MLE behaves similarly to the James--Stein benchmark, while NOMAD-CMLE still tends to provide modest gains. As the signal magnitude increases, the advantage of the NOMAD procedures becomes more pronounced.

\subsection{Linear regression simulations}\label{sec:sims_regression}

We generate data from the linear model in \eqref{eq:linreg_model} with $(n,p)=(1000,500)$ and noise variance $\sigma^2=1$. The design matrix $\pmb{X}$ is standardized and generated so that the eigenvalues of its correlation matrix are sampled from $U(0,10)$. The regression coefficient vector $\pmb{\beta}$ follows the same spike-and-slab model used in the normal mean experiments, with sparsity varied over a grid from $0$ to $1$. We vary the overall signal level so that the population coefficient of determination satisfies
$R^2\in\{0.25,0.50,0.75\}$,
which correspond to weak, moderate, and strong regression signal regimes.

\begin{figure}[h!]
    \centering
    \begin{subfigure}[b]{0.32\textwidth}
        \centering
        \includegraphics[width=\textwidth]{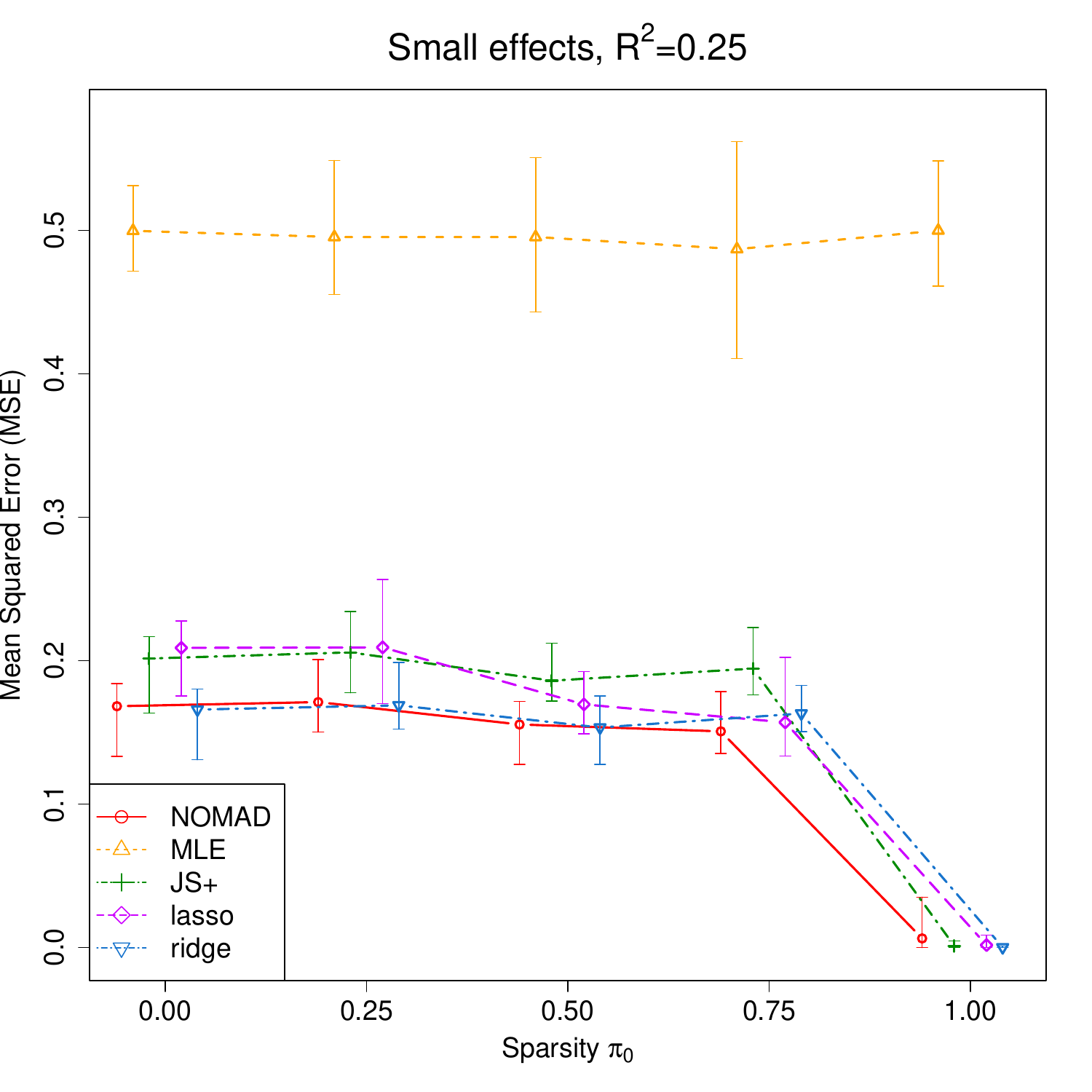}
    \end{subfigure}
    \begin{subfigure}[b]{0.32\textwidth}
        \centering
        \includegraphics[width=\textwidth]{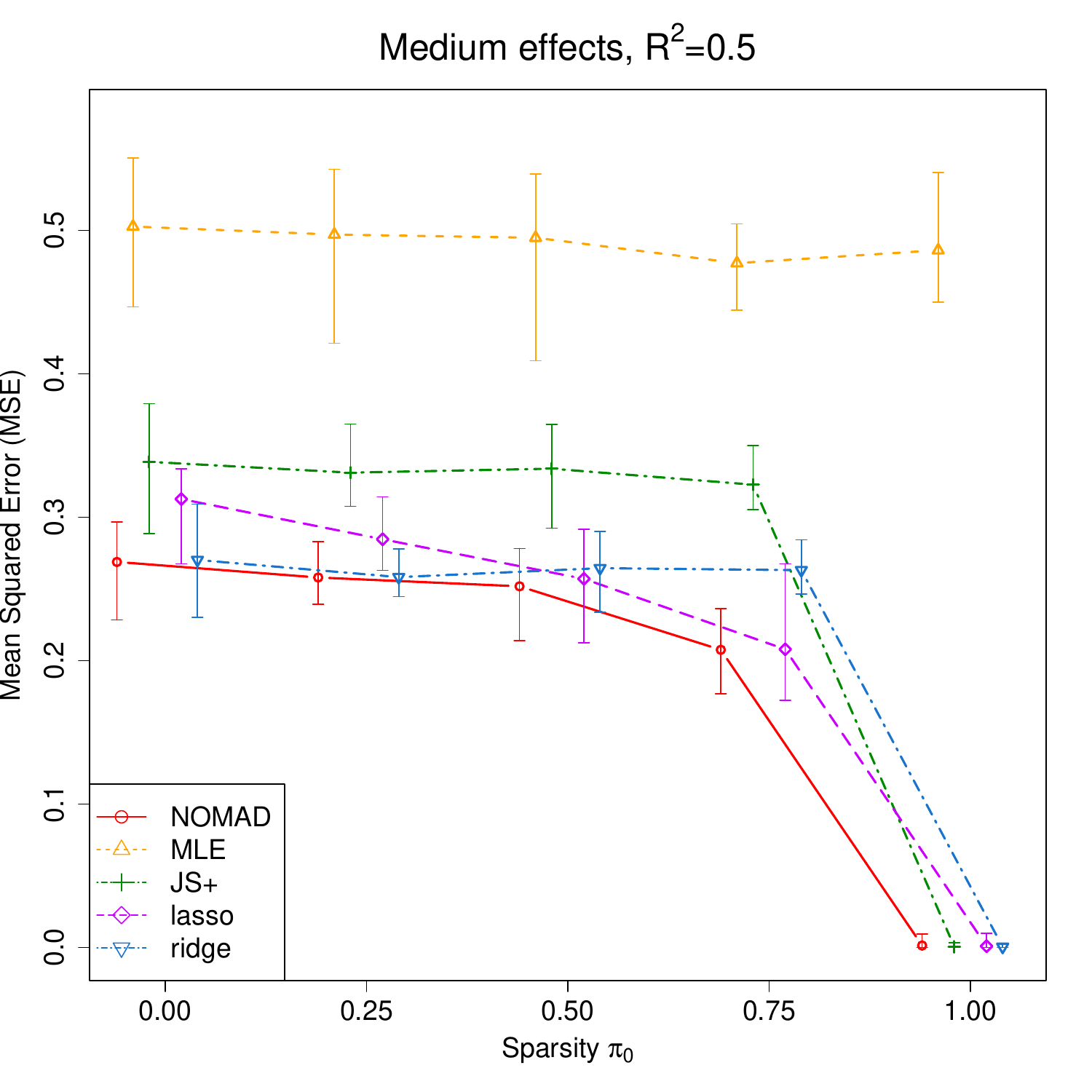}
    \end{subfigure}
    \begin{subfigure}[b]{0.32\textwidth}
        \centering
        \includegraphics[width=\textwidth]{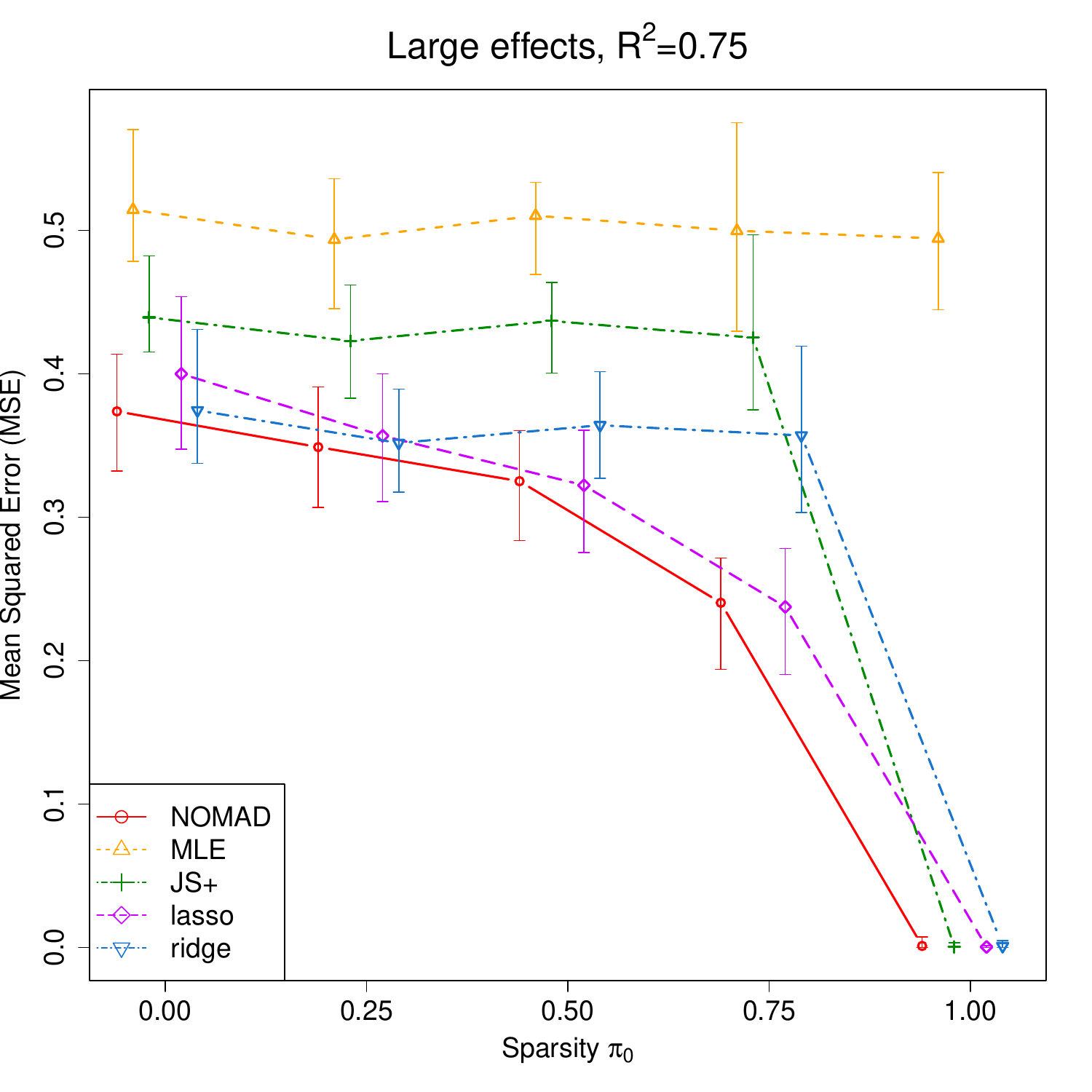}
    \end{subfigure}
    \caption{Prediction error in the linear regression simulations for weak, moderate, and strong signal regimes.}
    \label{fig:lr_mse}
\end{figure}

We compare:
(a) least squares;
(b) whitening-based positive-part James--Stein;
(c) lasso with cross-validated tuning;
(d) ridge with cross-validated tuning;
and (e) the NOMAD regression estimator,  
and report the prediction error
$
\frac{1}{n}\|\pmb{X}\hat{\pmb{\beta}}-\pmb{X}\pmb{\beta}\|_2^2
$
in Figure~\ref{fig:lr_mse}.

Least squares has the largest prediction error throughout. Ridge regression is consistently more stable than the James--Stein benchmark, and the gap widens as signal strength increases. Lasso performs especially well in sparse regimes and improves on ridge when the number of nonzero coefficients is small. NOMAD performs best overall across the range of sparsity and signal levels considered. In sparse regimes, its performance is close to that of lasso, while in denser regimes it behaves more like ridge. This pattern is consistent with the penalized-regression interpretation developed in Section~\ref{sec:regression}, where NOMAD induces a data-adaptive penalty that interpolates between shrinkage-dominant and thresholding-dominant behavior.





\section{Discussion}\label{sec:discussion}

This paper develops a unified framework for shrinkage-thresholding estimation based on approximate risk minimization. Starting from a general functional class that accommodates both shrinkage-dominant and thresholding-dominant behavior, the framework selects an estimator by minimizing a data-dependent approximation to quadratic risk. The resulting procedure, NOMAD, adapts to the empirical marginal regime of the sufficient statistics and provides a single construction that connects classical shrinkage rules with sparse thresholding rules.

Approximate risk minimization offers a common objective for methods that are often motivated separately. Shrinkage rules such as James--Stein are typically associated with dense regimes and global risk improvement, whereas thresholding rules such as soft-thresholding and lasso are associated with sparsity and variable selection. In the present framework, these behaviors arise as different regions of the same admissible class, and the score structure of the empirical marginal distribution determines the balance between them. This perspective clarifies why shrinkage is favored in dense regimes, why thresholding is favored in sparse regimes, and how intermediate behavior can be selected without imposing a regime-specific rule in advance.

The canonical normal mean model provides the cleanest setting for defining the estimator class, deriving the risk functional, and constructing the approximate risk criterion. Extending the framework to correlated multivariate normal mean estimation shows that the same principle can be combined with covariance or precision structure, which is essential for applications such as regression. In particular, the conditional construction illustrates how one can retain coordinatewise shrinkage-thresholding behavior while incorporating dependence through a stable normalization. The linear regression formulation further shows that the method admits an equivalent penalized representation with a data-adaptive penalty induced by the shrinkage-thresholding map.

The approximate risk criterion is the operational center of the method because it replaces infeasible oracle-risk optimization with a computable objective. Its construction highlights the importance of marginal score estimation and of distributional shape in determining the selected shrinkage-thresholding profile. The present theory establishes optimizer characterization and sieve-based stability in the canonical model, and consistency results in both canonical and correlated settings. At the same time, a more complete understanding of approximation error, convergence rates, and finite-sample behavior remains an important direction for future work.

Several limitations remain. First, the correlated extension requires covariance or precision input, and error in these quantities can affect both stability and performance, especially in high-dimensional settings \cite{cai2011clime}. Second, computation is driven by score estimation and by the optimization or fixed-point steps needed to determine the shrinkage-thresholding profile. Third, the Gaussian assumption is central to the normal mean reductions and Stein-type risk identities used in the analysis, so substantial departures from Gaussianity may reduce the accuracy of the approximate-risk construction.

Several extensions also appear natural. One direction is to adapt the framework beyond Gaussian models, for example through surrogate risk criteria for generalized linear models or robust loss-based formulations under non-Gaussian noise \cite{wang2012penalized,audibert2011robust,avella2017influence}. Another is to study uncertainty in covariance or precision estimation within the correlated construction \cite{cai2011clime}. It is also natural to investigate variable selection properties more systematically, including support recovery, false discoveries, and post-selection inference in sparse regimes \cite{wasserman2009high,su2017false,lee2016exact}.

Overall, the main contribution of this work is a unified approximate-risk framework for selecting shrinkage-thresholding rules in normal mean estimation problem. The framework provides a principled link among shrinkage, thresholding, and penalized regression, and it offers a flexible foundation for further methodological and theoretical development.

\begin{supplement}
\stitle{Supplement to "Approximate Risk Minimization over Shrinkage–Thresholding Rules in Normal Mean Estimation"}
\sdescription{Proof of Lemmas, Propositions and Theorems.}
\end{supplement}


\bibliographystyle{imsart-number} 
\bibliography{NOMAD}       




\end{document}